\documentclass[MSNbibl,number,citesort,seceqn,dvips]{arxbj}
\usepackage{upgreek}

\aid{0}
\volume{20}
\issue{3}
\pubyear{2014}
\firstpage{1059}
\lastpage{1096}
\doi{10.3150/13-BEJ515} 

\makeatletter

\newcommand{\RMi}{\mathrm{i}}

\renewcommand{\mid}{|}

\newcommand{\RMo}{\mathrm{o}}

\newcommand{\RMe}{\mathrm{e}}

\newcommand{\mrmd}{\,\mathrm{d}}
\newcommand{\rrvert}{\vert}
\newcommand{\llvert}{\vert}

\newtheorem{lem}{Lemma}[section]
\newtheorem{theo}{Theorem}[section]
\newtheorem{prop}{Proposition}[section]
\newtheorem{coro}{Corollary}[section]

\newremark{rem}{Remark}[section]

\newcommand{\cal}{\mathcal}

\newcommand{\xrightarrow}[2]{{}\mathop{\hbox to 1cm{\rightarrowfill}}^{#2}_{#1}{}}

\newcommand{\ind}{1}

\newcommand{\opt}{\mathrm{opt}}
\newcommand{\Obs}{\mathcal{O}}

\makeatother

\begin{document}
\begin{frontmatter}

\title{Asymptotic lower bounds in estimating jumps}
\runtitle{Asymptotic lower bounds in estimating jumps}

\begin{aug}
\author[1]{\inits{E.}\fnms{Emmanuelle} \snm{Cl\'ement}\corref{}\thanksref{1}\ead[label=e1]{emmanuelle.clement@univ-mlv.fr}},
\author[2]{\inits{S.}\fnms{Sylvain} \snm{Delattre}\thanksref{2}\ead[label=e2]{sylvain.delattre@math.jussieu.fr}} \and
\author[3]{\inits{A.}\fnms{Arnaud} \snm{Gloter}\thanksref{3}\ead[label=e3]{arnaud.gloter@univ-evry.fr}\ead[label=u1,url]{http://www.foo.com}}
\runauthor{E. Cl\'ement, S. Delattre and A. Gloter} 
\address[1]{Universit\'e Paris-Est, LAMA (UMR 8050), UPEMLV, CNRS,
F-77454, Marne-la-Vall\'ee, France. \printead{e1}}
\address[2]{Laboratoire de Probabilit\'es et Mod\`eles Al\'eatoires,
UMR 7599,
Universit\'e Paris Diderot, 75251 Paris Cedex 05, France. \printead{e2}}
\address[3]{Universit\'e d'\'Evry Val d'Essonne, D\'epartement de
Math\'ematiques,
91025 \'Evry Cedex, France.\\ \printead{e3}}
\end{aug}

\received{\smonth{2} \syear{2012}}
\revised{\smonth{9} \syear{2012}}

%
\begin{abstract}
We study the problem of the efficient estimation of the jumps for
stochastic processes. We assume that the stochastic jump process
$(X_t)_{t \in[0,1]}$ is observed discretely, with a sampling step of
size $1/n$. In the spirit of Hajek's convolution theorem, we show some
lower bounds for the estimation error of the sequence of the jumps
$(\Delta X_{T_k})_k$. As an intermediate result, we prove a LAMN
property, with rate $\sqrt{n}$, when the marks of the underlying jump
component are deterministic. We deduce then a convolution theorem, with
an explicit asymptotic minimal variance, in the case where the marks of
the jump component are random. To prove that this lower bound is
optimal, we show that a threshold estimator of the sequence of jumps
$(\Delta X_{T_k})_k$ based on the discrete observations, reaches the
minimal variance of the previous convolution theorem.
\end{abstract}

%

\pdfkeywords{convolution theorem,
Ito process, LAMN property}

\end{frontmatter}

\section{Introduction}

The statistical study of stochastic processes with jumps, from high
frequency data, has been the subject of many recent works.
A major issue is to determine if the jump part is relevant to model
the observed phenomenon. Especially, for modelling of asset prices, the
assessment of the part
due to the jumps in the price is an important question. This has been
addressed in several works, either by considering multi-power
variations \cite{BarShe06,BarSheWin06,HuaTau06}
or by truncation methods (see \cite{Mancini04,Mancini09}).
Another issue is to test statistically if the stochastic process
has continuous paths. The question has been addressed in many works
(see \cite{AitSahalia02,AitJia09,AitJac09AOSa})
and is crucial
to the hedging of
options.
A clearly related question is to determine the degree of activity of
the jump component of the process. Estimators of the Blumenthal--Getoor
index of the L\'evy measure of the process are proposed in several
papers \cite{AitJac09AOSb,ConMan11,Woerner06}.

In that context, the main statistical difficulty comes from the fact
that one observes a discrete sampling of the process, and consequently,
the exact values of the jumps are unobserved.
As a matter of fact, a lot of statistical procedures rely on
the estimation of a functional of the jumps.
In \cite{Jacod08SPA}, Jacod considers the estimation, from a high
frequency sampling, of the functional of the jumps
\[
\sum_{0\le s \le1, \Delta X_s \neq0} f(\Delta X_s) =\sum
_{k} f(\Delta X_{T_k})
\]
for a smooth function $f$ vanishing at zero (see Theorems 2.11 and 2.12
in \cite{Jacod08SPA} for precise assumptions).
In particular, he studies the difference between the unobserved
quantity $\sum_{0\le s \le1} f(\Delta X_s)$ and the observed one
$\sum_{i=0}^{n-1} f(X_{({i+1})/{n}}-X_{{i}/{n}})$. When $X$ is a
semi-martingale, it
is shown that the difference between the two quantities goes to zero
with rate $\sqrt{n}$. Rescaled by this rate, the difference is
asymptotically distributed as
%
\begin{equation}
\label{ElawJacod} \sum_{k}
f'(\Delta X_{T_k})\bigl[\sigma_{T_k-}
\sqrt{U_k} N_k^-+\sigma_{T_k}
\sqrt{1-U_k} N_k^+\bigr],
\end{equation}
where the variables $U_k$ are uniform variables on $[0,1]$ and $N_k^-$,
$N_k^+$ are standard Gaussian variables. The quantity $\sigma_{T_k-}$
(resp., $\sigma_{T_k}$) is the local volatility of the semi martingale
$X$ before (resp., after) the jump at time $T_k$.
This result
serves as the basis for studying the statistical procedures developed
in \cite{AitJac09AOSa,JacTod10}.

However, the problem of the efficiency of these methods seems to have
never been addressed.
Motivated by these facts,
we discuss, in this paper, the notion of efficiency to estimate the jumps
from the discrete sampling $(X_{i/n})_{0 \le i \le n}$.

Let us stress, that the meaning of efficiency is not straightforward
here. Indeed, we are not dealing with a standard parametric statistical
problem, and it is not clear which quantity can stand for the Fisher's
information.
In this paper, we restrict ourself to processes $X$ solutions of
\[
X_t=x_0+ \int_0^t
b(s,X_s) \mrmd  s + \int_0^t
a(s,X_s) \mrmd  W_s + \sum_{T_k \le t}
c(X_{T_{k}-}, \Lambda_k),
\]
where we assume that the number of jumps on $[0,1]$, denoted by $K$, is
finite. We
note
$J=(\Delta X_{T_1},\ldots,\Delta X_{T_K})$ the vector of jumps,
and $\Lambda=(\Lambda_1,\ldots,\Lambda_K)$ the random marks.
The notion of efficiency will be stated in this context
as a convolution result in Theorem \ref{Toptimaliteconv}.
More precisely, we prove that for any estimator $\widetilde{J}{}^{n}$
such that
the error $\sqrt{n}(\widetilde{J}{}^{n}-J)$ converges in law to some
variable $Z$, the law of $Z$ is a convolution between the law of
the vector
\[
\bigl[a(T_k,X_{T_k-}) \sqrt{U_k}
N_k^-+ a(T_k,X_{T_k}) \sqrt{1-U_k}
N_k^+\bigr]_{k=1,\ldots,K}
\]
and some other law. Contrarily to the standard convolution theorem, we
do not need the usual regularity assumption on the estimator. The
explanation is that we are not estimating a deterministic (unknown)
parameter, but we estimate some random (unobserved) variable $J$.

The proof of this convolution result relies on the study of a
preliminary parametric model: we consider the parametric model where
the values of the marks $\Lambda$ are considered as an unknown
deterministic parameter $\lambda\in\mathbb{R}^K$. The resulting model
is a stochastic differential equation with jumps, whose coefficients
depend on this parameter
$\lambda$. We establish then in Theorem
\ref{Th-LAMN}, that this statistical experiment satisfies the LAMN
property, with rate $\sqrt{n}$ and some explicit Fisher's information
matrix $I(\lambda)$.

By Hajek's theorem, it is well known that the LAMN property implies a
convolution theorem for any regular estimator of the parameter $\lambda
$ (see \cite{Jeganathan82,vandervaart}).
However, our context differs from the usual Hajek's convolution theorem
on at least two points.
First, the parameter $\lambda$ is randomized and second the target of
the estimator $J=(c ( X_{T_k-},\Lambda_k))_k$ depends both on the
randomized parameter and on some unobserved quantities $X_{T_k-}$. As a
result, the connection between the minimal law of the convolution
theorem and the Fisher's information of the parametric model is not
straightforward. The proof of the convolution theorem, when $c (
X_{T_k-},\Lambda_k)=c(\Lambda_k)$ does not depend on $X_{T_k-}$, is
simpler and is given in Theorem \ref{TConvosimple}.

Remark that it is certainly possible to state a general result about
the connection
between the LAMN property and convolution theorems for the estimation
of unobserved random quantities.
The proof of the Proposition \ref{PindepW} is a step in this direction.
However, giving such general results is beyond the scope of the paper.

The outline of the paper is as follows.
In Section~\ref{SOptim}, we state a convolution theorem, which
establishes an asymptotic lower bound
for the asymptotic error of any estimator of the jumps.
The LAMN property is enounced in Section~\ref{SLAMN}.
In Section~\ref{SEstimresult}, we
show that the threshold estimator, introduced by Mancini (see \cite{Mancini04,Mancini09}),
reaches the lower bound of Theorem
\ref{Toptimaliteconv}. This proves that this lower bound is optimal.
The proofs of these results are postponed to the Section~\ref{SProofs}.

\section{Convolution theorem}\label{SOptim}
\subsection{Notation}\label{SNotation}
Consider $(X_t)_{t \in[0,1]}$ an adapted c.\`a.d.l.\`a.g., one
dimensional, stochastic process defined on some filtered probability space
$(\Omega, \mathcal{F}, (\mathcal{F}_t)_{t \in[0,1]}, \mathbb{P})$. We
assume that the sample paths of $X$ almost surely admit a finite number
of jumps. We denote by $K$ the random number of jumps on $[0,1]$ and
$0<T_1 <\cdots<T_K<1$ the instants of these jumps. We assume that the
process $X$ is a solution of the stochastic differential equation with jumps
%
\begin{equation}
\label{Emodelbornesup} X_t=x_0+
\int_0^t b(s,X_s) \mrmd  s + \int
_0^t a(s,X_s) \mrmd
W_s + \sum_{T_k \le t} c(X_{T_{k}-},
\Lambda_k),
\end{equation}
where $W$ is a standard $(\mathcal{F}_t)_{t}$ Brownian motion.
The vector of marks $(\Lambda_k)_k$ is random. The Brownian motion, the
jump times and the marks are independent.

We will note $J=(J_k)_{k \geq1}$ the sequence of the jumps of the
process, defined by $J_k =c(X_{T_k-},\Lambda_k)=\Delta X_{T_k} $, for
$1 \leq k\leq K$ and $J_k=0$, for $k>K$.

Remark that if $T_k-T_{k-1}$ is exponentially distributed, the jumps
times are arrival times of a Poisson process.
Then, if the marks $(\Lambda_k)_{k}$ are i.i.d. variables, the process
$\sum_{T_k \le t} \Lambda_k$ is a compound Poisson process. In this
particular case, the equation
(\ref{Emodelbornesup}) becomes a standard SDE with jumps based on
a random Poisson measure with finite intensity.

It is convenient to assume that the process is realized on the
canonical product space of the Brownian part and the jumps parts
$\Omega
=\Omega^1 \times\Omega^2$,
$\mathbb{P}=\mathbb{P}^{1}\otimes\mathbb{P}^{2}$. More precisely, we
note $(\Omega^1, \mathcal{F}^{1}, \mathbb{P}^1)=(\mathcal{C}([0,1]),
\mathcal{B}, \mathbb{W})$, the space of continuous functions endowed
with the Wiener measure on the Borelian sigma-field and
$(\mathcal{F}^{1}_t)_{t \in[0,1]}$ the filtration generated by the
canonical process. We introduce $(\Omega^2, \mathcal{F}^{2},
\mathbb{P}^2)=(\mathbb{R}^{\mathbb{N}}\times
\mathbb{R}^{\mathbb{N}},\mathcal{B}(\mathbb{R})^{\otimes\mathbb
{N}}\otimes
\mathcal{B}(\mathbb{R})^{\otimes\mathbb{N}},\mathbb{P}^2)$, where
$\mathbb{P}^2$ is the law of two independent sequences of random
variables $(T_k)_{k \ge1}$, $(\Lambda_k)_{k \ge1}$. We assume that,
$\mathbb{P}^2$-almost surely, the sequence $(T_k)_{k \ge1}$ is
nondecreasing and such $K= \operatorname{Card}\{k, T_k \leq1\}$ is finite.
Then, $((W_t)_{t \in[0,1]},(T_k)_{k \ge1},(\Lambda_k)_{k \ge1})$ are
the canonical variables on $\Omega$. We assume that $(\mathcal{F}_t)_t$
is the right continuous, completed, filtration based on
$(\mathcal{F}^{1}_t \times\mathcal{F}^{2})_t$ and
$\mathcal{F}=\mathcal{F}_1$.

In order to describe the asymptotic law of any estimator of the jumps,
we need some additional notation. Following \cite{Jacod08SPA}, we
introduce an extension of our initial probability space. We consider an
auxiliary probability space $(\Omega', \mathcal{F'}, \mathbb{P'})$
which contains
$U=(U_k)_{k \ge1}$ a sequence of independent variables with uniform
law on $[0,1]$, and
$N^-=(N^{-}_k)_{k \ge1}$, $N^+=(N^{+}_k)_{k \ge1}$ two sequences of
independent variables with standard Gaussian law. All these
variables\vspace*{1pt}
are mutually independent.
We extend the initial probability space by setting $\widetilde{\Omega
}=\Omega\times\Omega'$,
$\widetilde{\mathcal{F}}=\mathcal{F}\otimes\mathcal{F'}$,
$\widetilde
{\mathbb{P}}= \mathbb{P}
\otimes\mathbb{P}'$,
$\widetilde{\mathcal{F}}_t=\mathcal{F}_t \otimes\mathcal{F}'$.

\subsection{Main result}

We need some more assumptions on the process. Especially, to avoid
cumbersome notation we will first assume in the next subsection that
the number of jumps is deterministic. We will show in Section~\ref
{SsRandom} that this is not a real restriction, since
we can reformulate our result by conditioning on the number of jumps $K$.

\subsubsection{Deterministic number of jumps}

Since the number of jumps $K$ is deterministic, the probability space
$\Omega$ introduced in Section~\ref{SNotation} is simplified
accordingly:
$\Omega=\Omega^1 \times\Omega^2$, $\Omega^1=\mathcal{C}([0,1])$ and
$\Omega^2=\mathbb{R}^K
\times\mathbb{R}^K$. The space $\widetilde{\Omega}=\Omega\times
\Omega
'$ with $\Omega'=\mathbb{R}^{3K}$ extends the initial space with the
sequences $N^-=(N_k^-)_{1\le k\le K}$, $N^+=(N_k^+)_{1\le k\le K}$,
$U=(U_k)_{1\le k\le K}$.

\renewcommand\thelonglist{H\arabic{longlist}}
\renewcommand\labellonglist{\thelonglist}
\begin{longlist}
\setcounter{longlist}{-1}
\item\label{HA0lower}
(\textit{Law of the
jump times}). The number of jumps $K$ is deterministic and the law of
$T=(T_1,\ldots, T_K)$ is absolutely continuous with respect to the
Lebesgue measure. We note $f_T$ its density.

\item\label{HA1lower} (\textit{Smoothness assumption}). The functions $(t,x)
\mapsto a(t,x)$ and $(t,x) \mapsto b(t,x)$ are $\mathcal{C}^{1,2}$ on
$[0,1] \times\mathbb{R}$. We note $a'$ and $b'$ their derivatives with
respect to $x$ and
we assume that $a'$ and $b'$ are $\mathcal{C}^{1,2}$ on $[0,1] \times
\mathbb{R}$. Moreover, the functions
$a$, $b$, and their derivatives are uniformly bounded.

The function $(x,\theta) \mapsto c(x,\theta)$ is $\mathcal{C}^{2,1}$ on
$\mathbb{R} \times\mathbb{R}$, with bounded derivatives. We note $c'$
its derivative with respect to $x$ and $\dot{c}$ its derivative with
respect to $\theta$. We assume moreover that $\dot{c}$
is $\mathcal{C}^{1,1}$ with bounded derivatives.

\item\label{HA2lower} (\textit{Non-degeneracy assumption}). We assume that
there exist two constants $\underline{a}$ and $\overline{a}$
such that
\begin{eqnarray*}
&\forall(t,x) \in[0,1] \times\mathbb{R}\qquad 0 < \underline{a} \leq a(t,x)
\leq
\overline{a};&
\\
&\forall(x,\theta) \in\mathbb{R} \times\mathbb{R}\qquad \bigl|1+c'(x,
\theta) \bigr|\geq\underline{a}.&
\end{eqnarray*}

\item\label{HA3lower} (``\textit{Randomness}'' \textit{of the jump
sizes}). The law of $\Lambda=(\Lambda_1,\ldots,\Lambda_K)$ is absolutely
continuous with respect to the Lebesgue measure and we note $f_\Lambda$
its density. We assume also
\[
\forall(x,\theta) \in\mathbb{R} \times\mathbb{R}\qquad \dot{c}(x,\theta)\neq0.
\]
\end{longlist}

Let us comment on these assumptions. First, the assumption that the
vector of jump times admits a density, hypothesis \ref{HA0lower}, is
crucial to prove the convergence in law of the fractional part of
$(nT_k)_k$ to the vector of uniform laws $(U_k)_k$.
In order to find a lower bound, we need to deal with a kind of regular
model, this explains
the assumption \ref{HA1lower}. Moreover, it is clear that if the
diffusion coefficient $a$ is equal to zero, one will expect a rate of
convergence for the estimation of the jumps faster than $\sqrt{n}$. In
that case, the LAMN property will not be satisfied with rate $\sqrt{n}$.
This clarifies why we assume a strictly positive lower bound on $a$.

Remark that the non-degeneracy of $\llvert 1+c'(x,\theta)\rrvert$ is
a standard
assumption which implies that the equation (\ref{Emodelbornesup})
admits a flow.
The assumption \ref{HA3lower} is more specifically related to our
statistical problem. We want to prove a lower bound for the estimation
of the random jump sizes.
Indeed, if these quantities do not exhibit enough randomness, it could
be possible to estimate them with a rate faster than $\sqrt{n}$. For
instance, the condition \ref{HA3lower} excludes that the jump sizes
do not depend on the underlying random marks.

We can now state our main result. We recall that
\[
J=(J_k)_{1\le k\le K}=\bigl(c(X_{T_k-},
\Lambda_k)\bigr)_{1\le k \le K}=(\Delta X_{T_k})_{1\le k \le K}
\in\mathbb{R}^K
\]
is the sequence of the jumps of the process.

We will call $(\widetilde{J}{}^{n})_{n \ge1}$ a sequence of estimators
if for each $n$, $\widetilde{J}{}^{n}\in\mathbb{R}^K$ is a measurable
function of the observations
$(X_{i/n})_{i=0,\ldots,n}$.

\begin{theo}\label{Toptimaliteconv}
Assume \textup{\ref{HA0lower}--\ref{HA3lower}}. Let $\widetilde{J}{}^{n}$ be
any sequence of estimators such that
%
\begin{equation}
\label{Epseudoregest} \sqrt{n}\bigl(
\widetilde{J}{}^{n}-J\bigr) \xrightarrow{\mathit{law}} {n \to\infty}
\widetilde{Z}
\end{equation}
for some variable $\widetilde{Z}$. Then, the law of $\widetilde{Z}$ is
necessarily a convolution:
%
\begin{equation}
\label{Econvolforme} \widetilde{Z}\stackrel{\mathit{law}} {=}
\bigl(\sqrt{U_k}a(T_k,X_{T_k-})
N_k^{-}+ \sqrt{1-U_k}a(T_k,X_{T_k})
N_k^{+} \bigr)_k + \widetilde{R},
\end{equation}
where conditionally on
$(T, \Lambda, (W_t)_{t\in[0,1]}, (U_k)_k)$, the random vector
$\widetilde{R}$ is independent of $(N_k^{-},N_k^{+})_k$.
\end{theo}

We will say that an estimator $\widetilde{J}{}^{n}$ of the jumps is
efficient if the asymptotic distribution
of $\sqrt{n}(\widetilde{J}{}^{n}-J)$ is equal in law to
$ (\sqrt{U_k}a(T_k,X_{T_k-}) N_k^{-}+ \sqrt{1-U_k}a(T_k,X_{T_k})
N_k^{+} )_k$ (which corresponds to $\widetilde{R}=0$).

It is well known that
in parametric models, the Hajek's convolution theorem usually requires
a regularity assumption on the estimator (see \cite{IbrHas81Book,vandervaart}). Here, our theorem does not require any assumption
on the estimator, apart its convergence with rate $\sqrt{n}$.
This comes from the fact that the target $J$ of the estimator is
random, yielding to some
additional regularity properties, compared with the usual parametric
setting (see a related situation in Jeganathan \cite{Jeganathan81}).
%
\begin{rem}\label{Ropti}
We can observe that
\[
\bigl(\sqrt{U_k}a(T_k,X_{T_k-})
N_k^{-}+ \sqrt{1-U_k}a(T_k,X_{T_k})
N_k^{+} \bigr)_k \stackrel{\mathrm{law}}
{=} \bigl(I^\opt\bigr)^{-1/2} N,
\]
where $I^\opt$ is the diagonal random matrix of size $K\times K$,
defined on the extended probability space $\widetilde{\Omega}$, with
diagonal entries:
%
\begin{equation}
\label{EdefI} I^\opt_{k}=
\bigl[U_k a(T_k,X_{T_k-})^2 +
(1-U_k)a(T_k,X_{T_k})^2
\bigr]^{-1}\qquad\mbox{for } k=1,\ldots,K.
\end{equation}
Conditionally on
$(T, \Lambda, (W_t)_{t\in[0,1]}, (U_k)_k)$, the vector $N$ is a
standard Gaussian vector on $\mathbb{R}^K$ and consequently $N$ is
independent of $I^\opt$.
\end{rem}

\begin{rem}\label{RIopt}
The Theorem \ref{Toptimaliteconv} states, in particular, that any
estimator of the jumps with rate $\sqrt{n}$ must have an asymptotic
conditional variance greater than $(I^\opt)^{-1}$.

Let us stress that if the rate of convergence is faster
than $\sqrt{n}$, then
(\ref{Epseudoregest}) is still true with $\widetilde{Z}=0$ and consequently
the Theorem \ref{Toptimaliteconv} proves that a convergence faster
than $\sqrt{n}$ is impossible.
\end{rem}

Now if instead estimating $J$, we estimate a function of the vector of
jumps, we can prove in a similar way the following
result. For the sake of shortness, we will omit the proof of the
following proposition.
%
\begin{prop}\label{Poptimaliteconv}
Assume \textup{\ref{HA0lower}--\ref{HA3lower}}.
Let $F$ be a ${\cal C}^1$ function from $\mathbb{R}^K$ to $\mathbb{R}$
and let $\widetilde{F}{}^{n}$ be any sequence of estimators of $F(J)$
such that
%
\begin{equation}
\label{Ecoropseudoregest} \sqrt{n}
\bigl(\widetilde{F}{}^{n}-F(J)\bigr) \xrightarrow{\mathit{law}} {n \to
\infty} \widetilde{Z}_F
\end{equation}
for some variable $\widetilde{Z}_F$. Then, the law of $\widetilde{Z}_F$
is necessarily a convolution:
%
\begin{equation}
\label{Ecoroconvolforme} \widetilde{Z}_F
\stackrel{\mathit{law}} {=} \sum_{k=1}^K
\frac{\partial F}{\partial x_k} (J) \bigl(\sqrt{U_k}a(T_k,X_{T_k-})
N_k^{-}+ \sqrt{1-U_k}a(T_k,X_{T_k})
N_k^{+} \bigr)+ \widetilde{R}_F,
\end{equation}
where, conditionally on
$(T, \Lambda, (W_t)_{t\in[0,1]}, (U_k)_k)$, the real random variable
$\widetilde{R}_F$ is independent of $(N_k^{-},N_k^{+})_k$.
\end{prop}
%
\begin{rem}\label{Roptijacod}
From the results of Jacod (Theorems 2.11 and 2.12 in
\cite{Jacod08SPA}), we deduce that the lower bound of Proposition \ref
{Poptimaliteconv} is optimal, and that the estimators of
\cite{Jacod08SPA} are efficient.
\end{rem}

\subsubsection{Random number of jumps} \label{SsRandom}
If the number of jumps is random, we need to modify some assumptions
accordingly.

\begin{longlist}[$\widetilde{\mathrm{H0}}$.]
\item[$\widetilde{\mathrm{H0}}$.]
We note $K=\operatorname{card} \{k
\mid T_k \in[0,1] \}$. Conditionally on $K$ the law of the vector of
jump times $T=(T_1, \ldots, T_K)$ admits a density.

\item[$\widetilde{\mathrm{H3}}$.] Conditionally on $K$, the law of
$(\Lambda_1,\ldots,\Lambda_K)$ is absolutely continuous with respect to
the Lebesgue measure. We assume also $
\forall(x,\theta) \in\mathbb{R} \times\mathbb{R}$ $\dot
{c}(x,\theta)\neq0$.
\end{longlist}

We can extend Theorem \ref{Toptimaliteconv}.
%
\begin{coro}
\label{CKrandom}
Assume $\widetilde{\mathrm{H0}}$, \textup{\ref{HA1lower},
\ref{HA2lower}}
and $\widetilde{\mathrm{H3}}$.
Let $\widetilde{J}{}^{n}$ be any sequence of estimators with values in
$\mathbb{R}^{\mathbb{N}}$ such that
\[
\sqrt{n}\bigl(\widetilde{J}{}^{n}-J\bigr) \xrightarrow{\mathit{law}} {n
\to\infty} \widetilde{Z}
\]
for some variable $\widetilde{Z}$. Then, the law of $\widetilde{Z}$
admits the decomposition:
\[
\widetilde{Z}\stackrel{\mathit{law}} {=} \bigl( \bigl[ \sqrt{U_k}a(T_k,X_{T_k-})
N_k^{-}+ \sqrt{1-U_k}a(T_k,X_{T_k})
N_k^{+} \bigr] 1_{\{ 1 \le k \le K\}} \bigr)_k
+ \widetilde{R},
\]
where conditionally on $(K, (T_k)_{1\le k \le K}, (\Lambda_k)_{1\le k
\le K}, (W_t)_{t\in[0,1]}, (U_k)_{1\le k \le K})$ the random vector
of the $K$ first components of $\widetilde{R}$ is independent of
$(N_k^{-},N_k^{+})_{1\le k \le K}$.
\end{coro}

In Section~\ref{SLAMN}, we consider a parametric model related to the
process (\ref{Emodelbornesup}), and enounce the associated LAMN
property. This is the key step before
proving Theorem \ref{Toptimaliteconv} and Corollary \ref{CKrandom}.
Remark that, directly considering the values of the jumps size as the
parameter, is not the right choice. The reason is that the jump sizes
are not independent of the Brownian motion $(W_t)_{t}$.
Instead, we prefer to consider the values of the marks $\Lambda$ as the
statistical parameter.

\section{LAMN property in an associated parametric model}\label{SLAMN}
We focus on the parametric model where the values of the marks $\Lambda
$ are considered as the unknown (deterministic) parameters, and $K$ is
deterministic. This is the crucial step before proving our convolution theorem.

More precisely, our aim is to obtain the LAMN property for the
parametric model
%
\begin{equation}
\label{Eq-par} X_t^{\lambda}=x_0+ \int
_0^t b\bigl(s,X_s^{\lambda}
\bigr) \mrmd  s + \int_0^t a\bigl(s,X_s^{\lambda}
\bigr) \mrmd  W_s + \sum_{k=1}^K
c\bigl(X_{T_{k}-}^{\lambda}, \lambda_k\bigr)
1_{t \geq T_k},
\end{equation}
where the parameter
$\lambda=(\lambda_1,\ldots, \lambda_K) \in\mathbb{R}^K$. We note
$T=(T_1,\ldots, T_K)$ the vector
of jump times such that $0<T_1< \cdots<T_K<1$.
Let us remark that, under the assumption \ref{HA0lower}, the
solutions of (\ref{Eq-par}) might be defined on the probability space
$\Omega^1 \times\mathbb{R}^K$ endowed with the product of the Wiener
measure and the law of the jumps times. But, to avoid new notation, we
can assume that, for all $\lambda\in\mathbb{R}^K$, the process
$(X^\lambda_t)_{t \in[0,1]}$ is defined on the space $(\Omega,
\mathcal
{F}, \mathbb{P} )$ of Section~\ref{SOptim}.

In this model, we assume that
we observe both the regular discretization $(X^{\lambda}_{i/n})_{1
\leq i \leq n}$ of the process solution of (\ref{Eq-par}) on the time
interval $[0,1]$ and the jump times vector $T$.
The observation of $T$ leads to a more tractable computation of the
likelihood. This is not restrictive to add some observations to the
statistical experiment,
since our aim is to derive an asymptotic lower bound.
Under \ref{HA0lower} and \ref{HA1lower}, the law of the
observations $(T,(X^{\lambda}_{i/n})_{1 \leq i \leq n})$ admits a density
$\mathbf{p}^{n,\lambda}$.
We note $\mathbf{p}^{n,\lambda,T}$ the density of $(X^{\lambda
}_{i/n})_{1 \leq i \leq n}$ conditionally on $T$.
For $h=(h_1,\ldots, h_K) \in\mathbb{R}^K$ we introduce the
log-likelihood ratio:
%
\begin{equation}\label{Eq-Zn}
Z_n(\lambda, \lambda+h/\sqrt{n}, T, x_1,\ldots,
x_n)=\log\frac
{\mathbf{p}^{n, \lambda+h/\sqrt{n}}}{ \mathbf{p}^{n, \lambda} }(T,x_1,\ldots,
x_n).
\end{equation}

\begin{theo}\label{Th-LAMN}
Assume \textup{\ref{HA0lower}, \ref{HA1lower}} and \textup{\ref{HA2lower}}. Then,
the statistical experiment $(\mathbf{p}^{n,\lambda})_{\lambda\in
\mathbb{R}^K}$ satisfies a LAMN property.
For $\lambda\in\mathbb{R}^K$, $h \in\mathbb{R}^K$ we have:
%
\begin{eqnarray}\label{Eq-LAMN}
&&
Z_n\bigl(\lambda, \lambda+h/\sqrt{n},
T, X^{\lambda}_{1/n},\ldots, X^{\lambda}_{1}\bigr)\nonumber\\[-8pt]\\[-8pt]
&&\quad= \sum_{k=1}^K
h_k I_n(\lambda)_k^{1/2}
N_n(\lambda)_k - \frac{1}{2} \sum
_{k=1}^K h_k^2
I_n(\lambda)_k +\RMo_{\mathbf
{p}^{n,\lambda}}(1),\nonumber
\end{eqnarray}
where $I_n(\lambda)$ is a diagonal random matrix and $N_n(\lambda)$ are
random vectors in $\mathbb{R}^K$ such that
\[
\bigl(I_n(\lambda),N_n(\lambda)\bigr) \xrightarrow{
\mathit{law}} {n \to\infty} \bigl(I(\lambda),N\bigr)
\]
with:
%
\begin{equation}
\label{EIpara} I(\lambda)_k=\frac{\dot{c}(X^{\lambda}_{T_k-},\lambda_k)^2}{
a^2(T_k,X^{\lambda}_{T_k-}) [1+c'(X^{\lambda}_{T_k-},\lambda_k)]^2 U_k+
a^2(T_k,X^{\lambda}_{T_k-}+c(X^{\lambda}_{T_k-},\lambda_k)
)(1-U_k)},\hspace*{3pt}
\end{equation}
where $U=(U_1,\ldots, U_K)$ is a vector of independent uniform laws on
$[0,1]$ such that $U$, $T$ and $(W_t)_{t \in[0,1]}$ are independent, and
conditionally on $(U,T,(W_t)_{t \in[0,1]})$,
$N$ is a standard Gaussian vector in $\mathbb{R}^K$.
\end{theo}
Actually, we can complete the statement of the theorem by giving
explicit expressions for $I_n(\lambda)$ and $N_n(\lambda)$:
\begin{eqnarray*}
I_n(\lambda)_k&=&\frac{ \dot{c}(X^\lambda_{{i_k}/{n}},\lambda
_k)^2}{nD^{n, \lambda_k,k}( X^\lambda_{{i_k}/{n}})},\\
N_n(\lambda)_k&=&\frac{\sqrt{n}(X^\lambda_{({i_k+1})/{n}}-X^\lambda_{{i_k}/{n}}-c(X^\lambda_{i_k/n},\lambda_k))} {
\sqrt{nD^{n, \lambda_k,k}( X^\lambda_{{i_k}/{n}}) }},
\\
D^{n,\lambda_k,k}\bigl( X^\lambda_{{i_k}/{n}}\bigr)&=&a^2
\biggl(\frac{i_k}{n}, X^\lambda_{{i_k}/{n}}\biggr)
\bigl(1+c'\bigl( X^\lambda_{{i_k}/{n}},\lambda
_k\bigr)\bigr)^2 \biggl(T_k-
\frac{i_k}{n}\biggr)
\\
&&{} + a^2\biggl(\frac{i_k}{n}, X^\lambda_{{i_k}/{n}}
+c\bigl( X^\lambda_{{i_k}/{n}},\lambda_k\bigr)\biggr)
\biggl(\frac{i_k+1}{n}-T_k\biggr),
\end{eqnarray*}
where $i_k$ is the integer part of $nT_k$.
%
\begin{rem}\label{RLAMN}
We remark that from a direct application of Hajek's theorem (see Van
der Vaart~\cite{vandervaart}, Corollary 9.9, page 132), any regular
estimator of $\lambda$ has an asymptotic conditional variance greater
than $I(\lambda)^{-1}$. Here, an estimator of $\lambda$ is a measurable
function of $(T,(X^{\lambda}_{i/n})_{1 \leq i \leq n})$, and so we
deduce that, a fortiori, any measurable function of $(X^{\lambda
}_{i/n})_{1 \leq i \leq n}$ satisfies the same asymptotic lower bound.
\end{rem}

\section{Efficient estimator of the jumps} \label{SEstimresult}

We use the notation of Section~\ref{SNotation} and since we just
propose an estimator of the jumps, we can weaken the assumptions of the
previous sections.

\renewcommand\thelonglist{A\arabic{longlist}}
\renewcommand\labellonglist{\thelonglist}
\begin{longlist}
\setcounter{longlist}{0}
\item\label{Hexistestcoef}
(\textit{Smoothness assumption}). The functions $a\dvtx
[0,1]\times\mathbb{R} \to\mathbb{R}$, $b\dvtx  [0,1]\times\mathbb{R} \to
\mathbb{R}$ and $c\dvtx  \mathbb{R}^2 \to\mathbb{R}$ are continuous.

\item\label{Hpositjumpestcoef} (\textit{Identifiability of the jumps}).
We have almost surely: $c(X_{T_k-},\Lambda_k) \neq0,\forall k \in\{
1,\ldots,K \}$.
\end{longlist}

This last condition ensures that the jump times of $X$ are exactly the
times $T_k$.

Recall that $J=(J_k)_{k \ge1}$ is the sequence of jumps of $X$ (on
$[0,1]$): we set
$J_k=\Delta X_{T_k}=X_{T_k}-X_{T_k-}$ for $k \le K$ and we define
$J_k=0$ for $k >K$.

We construct an estimator of $J$ following the threshold estimation
method proposed by Mancini \cite{Mancini04,Mancini09}.

Let $(u_n)_n$ be a sequence of positive numbers tending to $0$. We set
$\hat{i}_1^{n}= \inf\{0\le i \le n-1\mbox{:}\allowbreak\llvert X_{
({i+1})/{n}}-X_{i/n}\rrvert \ge u_n \}$ with the convention
$\inf
\varnothing= + \infty$. We recursively define for $k \ge2$,
%
\begin{equation}
\label{Edefihat} \hat{i}_k^{n}=
\inf\bigl\{ \hat{i}_{k-1}^{n}<i \le n-1\dvt\llvert
X_{({i+1})/{n}}-X_{i/n}\rrvert\ge u_n \bigr\}.
\end{equation}
We set $\widehat{K}_n=\sup\{ k \ge1\dvt\hat{i}^{n}_k<\infty\}$ the
number of increments of the jump diffusion exceeding the threshold $u_n$.
We then define for $k\ge1$,
%
\begin{equation}
\label{EdefJhat} \widehat{J}_k^{n}=
\cases{\displaystyle X_{{\hat{i}^{n}_k+1}/{n}}-X_{{\hat{i}^{n}_k}/{n}},
&\quad if
$k \le\widehat{K}_n$,
\cr
0, &\quad if $k > \widehat{K}_n$.}
\end{equation}
The sequence $(\widehat{J}{}^{n})_n$ is an estimator of the vector of
jumps $J$, and
$(\widehat{K}_n)_n$ estimates the number of jumps.

\begin{prop} \label{Pconsistence}
Let us assume $\widetilde{H}_0$, \textup{\ref{Hexistestcoef}, \ref
{Hpositjumpestcoef}} and $u_n \sim n^{-\varpi}$ with
$\varpi\in(0,1/2)$.
Then, we have almost surely,
\begin{eqnarray*}
&&\widehat{K}_n = K \qquad\mbox{for $n$ large enough}
\\
&&\qquad\mbox{if $k \le K$}\qquad \widehat{J}{}^{n}_k \xrightarrow{}
{ n \to\infty} J_k=\Delta X_{T_k},
\\
&&\qquad\mbox{if $k > K$}\qquad \widehat{J}{}^{n}_k = 0 \qquad\mbox{for
$n$ large enough}.
\end{eqnarray*}
\end{prop}

The consistency result concerning the estimator $\widehat{K}_n$ is a
special case of Mancini (\cite{Mancini09}, Theorem~1) and the jump
sizes $(J_k)_k$ were consistently estimated in Mancini \cite{Mancini04}
with exactly the same estimator but
when the observation time goes to infinity.

We now describe the asymptotic law of the error between $\widehat
{J}{}^{n}$ and $J$. Note that Theorem 3 in \cite{Mancini09} gives the
asymptotic distribution of the estimator of the sum of the jumps
assuming that the diffusion coefficient $a$ is independent of the
Brownian process $W$ and the jump part, this is the reason why the
uniform laws do not appear in the asymptotic law. The situation is
completely different here, since the diffusion coefficient $a$ depends
on the process $X$, and is more related
to Jacod's results (see \cite{Jacod08SPA}).

\begin{theo} \label{TTCLsaut}
Let us assume $\widetilde{H}_0$, \textup{\ref{Hexistestcoef}, \ref
{Hpositjumpestcoef}} and $u_n \sim n^{-\varpi}$ with
$\varpi\in(0,1/2)$.
Then $\sqrt{n}(\widehat{J}{}^{n}-J)$ converges in law to $Z=(Z_k)_{k
\ge
1}$ where
the limit can be described on the extended space $\widetilde{\Omega}$ by:
\begin{eqnarray*}
Z_k&=&\sqrt{U_k} a(T_k,X_{T_k-})
N_k^{-} + \sqrt{1-U_k}
a(T_k,X_{T_k}) N_k^{+}\qquad \mbox{for } k \le K,
\\
Z_k&=&0\qquad \mbox{for } k> K.
\end{eqnarray*}
Moreover\vspace*{1pt} the convergence is stable with respect to the sigma-field
$\mathcal{F}$.
Let us precise that, here, the convergence in law of the infinite
dimensional vector $\sqrt{n}(\widehat{J}{}^{n}-J)$ means the convergence
of any finite dimensional marginals.
\end{theo}

\begin{rem} \label{Rrandomvariance}
The Theorem \ref{TTCLsaut} shows that the error for the estimation of
the jump $\Delta X_{T_k}$ is asymptotically conditionally Gaussian and
that the estimator $\widehat{J}{}^{n}$ is efficient. In particular, the
conditional variance on $(T, \Lambda, K,(W_t)_{t \in[0,1]}, (U_k)_k)$
of the error is equal to the lower bound $(I^{\opt})^{-1}=U_k
a(T_k,X_{T_k-})^2 +(1-U_k)a(T_k,X_{T_k})^2$, and consequently this
lower bound is optimal.
\end{rem}

\section{Proof section}\label{SProofs}

We divide the proofs into three sections.

We first prove the LAMN property of the parametric model in Section~\ref
{SsLAMNpf}.
Then, the convolution result is established in Section~\ref{Ssconvolutionpf}.
Finally, the Section~\ref{Ssconsistencepf} is devoted to the proof
of the convergence and normality of the estimator $\widehat{J}{}^{n}$.

We first state a lemma which will be useful in the next sections.
%
\begin{lem}\label{L-tcltime}
Let $K_0 \in\mathbb{N}\setminus\{0\}$ and consider
$T=(T_1,\ldots,T_{K_0})$ a random variable on $[0,1]^{K_0}$ with density
$f_T$. For $k=1,\ldots,K_0$, we note $i_k=[nT_k]$ the integer part of
$nT_k$. Let $(W_t)_{t\in[0,1]}$ be a standard Brownian motion
independent of $T$.

Then, we have the convergence in law of the variables
\[
\biggl(T, \biggl(n\biggl(T_k-\frac{i_k}{n}\biggr)
\biggr)_{k}, \bigl(\sqrt{n}(W_{T_k}-W_{i_k/n})
\bigr)_k, \bigl(\sqrt{n}(W_{({i_k+1})/{n}}-W_{T_k})
\bigr)_k, (W_t)_{t \in[0,1]}\biggr)
\]
to
\[
\bigl(T, (U_k)_{k},\bigl(\sqrt{U_k}
N_k^-\bigr)_k, \bigl(\sqrt{1-U_k}
N_k^+\bigr)_k,(W_t)_{t
\in[0,1]}\bigr),
\]
where $U=(U_1,\ldots, U_{K_0})$ is a vector of independent uniform
laws on $[0,1]$, $N^-=(N_1^-,\ldots,\allowbreak N_{K_0}^-)$ and
$N^+=(N_1^+,\ldots, N_{K_0}^+)$ are independent standard Gaussian vectors
such that $T$, $U$, $N^-$, $N^+$ and $(W_t)_t$ are independent.
\end{lem}
\begin{pf}
The convergence of the vector
\[
\bigl(T, \bigl(\sqrt{n}(W_{T_k}-W_{{i_k}/{n}})
\bigr)_k, \bigl(\sqrt{n}(W_{(i_k+1)/{n}}-W_{T_k})
\bigr)_k, (W_t)_{t \in[0,1]}\bigr)
\]
is a direct consequence of Lemma 6.2 in \cite{JacPro98} (see also Lemma
5.8 in \cite{Jacod08SPA}) and following this proof (which is simpler
in our case), there is no difficulty to add the variables
$(n(T_k-\frac{i_k}{n}))_{k}$ in the vector.
\end{pf}

\subsection{LAMN property: Proof of Theorem \texorpdfstring{\protect\ref{Th-LAMN}}{3.1}}
\label{SsLAMNpf}

We use the framework of Section~\ref{SLAMN} and we
introduce some more notation. For $k=1,\ldots,K$, we note $i_k=[nT_k]$
the integer part of $nT_k$ and for $t \in[i_k/n, (i_k+1)/n]$, we note
$(X_t^{\theta,k})$ the process solution of the following jump-diffusion
equation with only one jump at time $T_k$:
%
\begin{equation}
\label{Eq-par-thetak} X_t^{\theta,k}=X_0+ \int
_{0}^t b\bigl(s,X_s^{\theta,k}
\bigr) \mrmd  s + \int_{0}^t a\bigl(s,X_s^{\theta,k}
\bigr) \mrmd  W_s + c\bigl(X_{T_{k}-}^{\theta,k}, \theta
\bigr)1_{t \geq T_k}.
\end{equation}
Under \ref{HA1lower} and \ref{HA2lower} and conditionally on $T$,
this process admits a strictly positive conditional density, which is
$\mathcal{C}^1$ with respect to $\theta$. We will note
$p^{\theta,T}(\frac{i_k}{n}, \frac{i_k+1}{n }, x,y)$ the density of
$X_{({i_k+1})/{n}}^{\theta,k}$ conditionally on $T$ and $X_{i_k/n}^{\theta,k}=x$ and $\dot{p}^{\theta,T}(\frac{i_k}{n},
\frac{i_k+1}{n }, x,y)$ its derivative with respect to $\theta$.

We observe that the log-likelihood ratio $Z_n$ only involves the
transition densities
of $X^\lambda$ on a time interval where a jump occurs. This transition is
$p^{\theta,T}(\frac{i_k}{n}, \frac{i_k+1}{n }, x,y)$ if there is
exactly one jump in the corresponding interval.
Then, one can easily see that the following decomposition holds for $Z_n$:
%
\begin{eqnarray}\label{Eq-Znk}
&&
Z_n\biggl(\lambda, \lambda+\frac{h}{\sqrt{n}}, T,
x_1,\ldots, x_n\biggr)\ind_{\mathcal{T}_n}(T)\nonumber \\
&&\quad=\sum
_{k=1}^K \ln\frac{p^{\lambda_k+{h_k}/{\sqrt
{n}},T}}{p^{\lambda
_k,T} } \biggl(
\frac{i_k}{n}, \frac{i_k+1}{n }, x_{i_k},x_{i_k+1}\biggr)
\ind_{\mathcal{T}_n}(T)
\\
&&\quad=\sum_{k=1}^K
\int_{\lambda_k}^{\lambda_k+
{h_k}/{\sqrt{n}}} \frac{\dot{p}^{\theta,T}}{p^{\theta,T} } \biggl(
\frac{i_k}{n}, \frac{i_k+1}{n }, x_{i_k},x_{i_k+1}\biggr)\mrmd
\theta\ind_{\mathcal{T}_n}(T),\nonumber
\end{eqnarray}
where $\ind_{\mathcal{T}_n}(T)$ is the indicator function that there is
at most one jump in each time interval $[i/n, (i+1)/n)$ for $i=0,\ldots,n-1$.

We have now to study the asymptotic behaviour of (\ref{Eq-Znk}). This
is divided into several lemmas. The Lemmas \ref{L-malliavin}--\ref
{L-reste} give an expansion for the
score function, with an uniform control in $\theta$. We deduce then an
explicit expansion for
$ \int_{\lambda_k}^{\lambda_k+{h_k}/{\sqrt{n}}}
\frac{\dot{p}^{\theta,T}}{p^{\theta,T} }
(\frac{i_k}{n}, \frac{i_k+1}{n },\break x_{i_k},x_{i_k+1})\mrmd  \theta$
in Lemma \ref{L-LAMN}, and conclude by passing through the limit in
Lemma \ref{L-tcl}.

We begin with
a representation of
$\frac{\dot{p}^{\theta,T}}{p^{\theta,T} }
(\frac{i_k}{n}, \frac{i_k+1}{n }, x,y)$ as a conditional expectation,
using Malliavin calculus. We refer to Nualart \cite{Nualart} for a
detailed presentation of Malliavin calculus.
The Malliavin calculus techniques to derive LAMN properties have been
introduced by Gobet \cite{Gobet01} in the case of
multi-dimensional diffusion processes and then used by Gloter and Gobet
\cite{GloGob08} for integrated diffusions.

In all what follows, we will denote by $C_p$ a constant (independent on
$n$, $k$ and $\theta$) which value may
change from line to line.

\begin{lem}\label{L-malliavin}
Assuming \textup{\ref{HA1lower}} and \textup{\ref{HA2lower}}, we have $\forall(x,y)
\in\mathbb{R}^2$:
\[
\frac{\dot{p}^{\theta,T}}{p^{\theta,T} } \biggl(\frac{i_k}{n}, \frac
{i_k+1}{n }, x,y
\biggr)=E^{x,T,k} \bigl(\delta\bigl(P^{n,
\theta,k}\bigr) |
X_{(i_k+1)/{n}}^{\theta,k}=y \bigr),
\]
where $E^{x,T,k}$ is the conditional expectation on $T$ and
$X^{\theta,k}_{{i_k}/{n}}=x$, $\delta$ is the Malliavin divergence
operator and $P^{n, \theta,k}$ is the process given on $[\frac{i_k}{n},
\frac{i_k+1}{n}]$ by
\[
P^{n, \theta,k}_s=\frac{ (Y_{T_k}^{\theta,k} Y_s^{\theta,k})^{-1}
(1+c'(X_{T_k-}^{\theta,k}, \theta)1_{s\leq T_k})
a(s, X_s^{\theta,k}) \dot{c}(X_{T_k-}^{\theta,k}, \theta)} {
\int_{i_k/n}^{(i_k+1)/n} (Y_u^{\theta,k})^{-2} a^2(u, X_u^{\theta,k})
(1+c'(X_{T_k-}^{\theta,k}, \theta)1_{u\leq T_k})^2 \mrmd u},
\]
where $(Y_t^{\theta,k})_t$ is the process solution of
%
\begin{equation}
\label{Eq-exp} Y_t^{\theta,k}=1+ \int_0^t
b'\bigl(s, X_s^{\theta,k}\bigr)Y_s^{\theta,k}\mrmd s
+ \int_0^t a'\bigl(s,
X_s^{\theta,k}\bigr)Y_s^{\theta,k}
\mrmd W_s.
\end{equation}
\end{lem}
We remark that under \ref{HA1lower}, the process $(Y_t^{\theta,k})_t$
and its inverse satisfy $\forall p \geq1$,
%
\begin{equation}
\label{Eq-BY} \Bigl(E\Bigl(\sup_{0\leq t\leq1} \bigl|
Y_t^{\theta,k}\bigr|^p\Bigr)\Bigr)^{1/p} \leq
C_p,\qquad \Bigl(E\Bigl(\sup_{0\leq t\leq1} \bigl|
Y_t^{\theta,k}\bigr|^{-p}\Bigr)\Bigr)^{1/p}
\leq C_p.
\end{equation}
\begin{pf}
The proof of Lemma \ref{L-malliavin} is based on Malliavin calculus on
the time interval $[i_k/n, (i_k+1)/n]$, conditionally on $T$ and
$(W_t)_{t \leq i_k/n}$. We first observe that under \ref{HA1lower}
and \ref{HA2lower}, the process $(X_t^{\theta,k})$ solution of
(\ref{Eq-par-thetak}) admits a derivative with respect to $\theta$
that we
will denote by $(\dot{X}_t^{\theta,k})$ (see, e.g., Kunita \cite
{Kunita} since this problem is similar to the derivative with respect
to the initial condition). Moreover
$(X_t^{\theta,k})$ and $(\dot{X}_t^{\theta,k})$ belong, respectively, to
the Malliavin spaces $\mathbb{D}^{2,p}$ and $\mathbb{D}^{1,p}$,
$\forall p \geq1$.
Now, let $\varphi$ be a smooth function with compact support, we have:
\[
\frac{\partial}{\partial\theta} E^{x,T,k} \varphi\bigl(X_{(i_k+1)/{n}}^{\theta,k}
\bigr) = E^{x,T,k} \varphi'\bigl(X_{(i_k+1)/{n}}^{\theta,k}
\bigr) \dot{X}_{(i_k+1)/{n}}^{\theta,k}.
\]
Using the integration by part formula (see Nualart \cite{Nualart},
Proposition 2.1.4, page 100), we can write
\[
E^{x,T,k} \varphi'\bigl(X_{(i_k+1)/{n}}^{\theta,k}
\bigr) \dot{X}_{(i_k+1)/{n}}^{\theta,k}= E^{x,T,k} \varphi
\bigl(X_{(i_k+1)/{n}}^{\theta,k}\bigr) H\bigl(X_{(i_k+1)/{n}}^{\theta,k},
\dot{X}_{(i_k+1)/{n}}^{\theta,k}\bigr),
\]
where the\vspace*{1pt} weight $H$ can be expressed in terms of the Malliavin
derivative of $X_{(i_k+1)/{n}}^{\theta,k}$, the inverse of its Malliavin
variance--covariance matrix and the divergence operator as follows:
\[
H\bigl(X_{(i_k+1)/{n}}^{\theta,k}, \dot{X}_{(i_k+1)/{n}}^{\theta,k} \bigr) =\delta\bigl(\dot{X}_{(i_k+1)/{n}}^{\theta,k} \gamma^{\theta,k} D X_{(i_k+1)/{n}}^{\theta,k}
\bigr),
\]
where
%
\begin{equation}
\label{Eq-gamma} \gamma^{\theta,k}= \biggl(\int_{i_k/n}^{(i_k+1)/n}
\bigl( D_u X_{(i_k+1)/{n}}^{\theta,k} \bigr)^2 \mrmd u
\biggr)^{-1}.
\end{equation}
On the other hand, from Lebesgue derivative theorem, we have:
\[
\frac{\partial}{\partial\theta} E^{x,T,k} \varphi\bigl(X_{(i_k+1)/{n}}^{\theta,k}
\bigr) =\int\varphi(y) \dot{p}^{\theta,T}\biggl(\frac{i_k}{n},
\frac{i_k+1}{n }, x,y\biggr)\mrmd y,
\]
this leads to the following representation
\begin{eqnarray*}
&&
\dot{p}^{\theta,T}\biggl(\frac{i_k}{n}, \frac{i_k+1}{n }, x,y
\biggr)\\
&&\quad=E^{x,T,k}\bigl(\delta\bigl(\dot{X}_{(i_k+1)/{n}}^{\theta,k}
\gamma^{\theta,k} D X_{(i_k+1)/{n}}^{\theta,k} \bigr)|
X_{(i_k+1)/{n}}^{\theta,k}=y\bigr)p^{\theta,T}\biggl(\frac{i_k}{n},
\frac{i_k+1}{n }, x,y\biggr).
\end{eqnarray*}
It remains to give a more tractable expression of $\dot{X}_{(i_k+1)/{n}}^{\theta,k} \gamma^{\theta,k} D X_{(i_k+1)/{n}}^{\theta,k} $. We first observe that:
\begin{eqnarray*}
\dot{X}_{(i_k+1)/{n}}^{\theta,k}&=&\dot{c}\bigl(X_{T_k-}^{\theta,k},
\theta\bigr) + \int_{T_k}^{(i_k+1)/n} b'
\bigl(u, X_u^{\theta,k}\bigr) \dot{X}_u^{\theta,k}
\mrmd u \\
&&{}+\int_{T_k}^{(i_k+1)/n} a'\bigl(u,
X_u^{\theta,k}\bigr) \dot{X}_u^{\theta,k}\mrmd W_u
\end{eqnarray*}
and consequently
%
\begin{equation}
\label{Eq-Xpoint} \dot{X}_{(i_k+1)/{n}}^{\theta,k}=Y_{(i_k+1)/{n}}^{\theta,k}
\bigl(Y_{T_k}^{\theta,k}\bigr)^{-1} \dot{c}
\bigl(X_{T_k-}^{\theta,k}, \theta\bigr),
\end{equation}
where $(Y_t^{\theta,k})$ is solution of (\ref{Eq-exp}). Turning to the
Malliavin derivative of $X_{(i_k+1)/{n}}^{\theta,k}$, we first
observe that $DX_{(i_k+1)/{n}}^{\theta,k} \in L^2([i_k/n,
(i_k+1)/n])$ and so we just have to
explicit $D_s X_{(i_k+1)/{n}}^{\theta,k}$ for $s \neq T_k$.
Assuming first that $T_k<s \leq(i_k+1)/n $,
we have for
$u \in[s, (i_k+1)/n] $:
\[
D_s X_{u}^{\theta,k}=a\bigl(s,
X_s^{\theta,k}\bigr) + \int_s^{u}
b'\bigl(v, X_v^{\theta,k}\bigr) D_s
X_v^{\theta,k} \mrmd v +\int_s^{u}
a'\bigl(v, X_v^{\theta,k}\bigr) D_s
X_v^{\theta,k}\mrmd W_v
\]
and then $D_s
X_{u}^{\theta,k}=Y_{u}^{\theta,k}(Y_{s}^{\theta,k})^{-1}a(s,
X_s^{\theta,k}) $.

Now, if $i_k/n \leq s<T_k$, we have for $u \geq s$
\begin{eqnarray*}
D_s X_{u}^{\theta,k}&=&a\bigl(s,
X_s^{\theta,k}\bigr) + c'\bigl(X_{T_k-}^{\theta,k},
\theta\bigr)D_s X_{T_k-}^{\theta,k} 1_{ u \geq T_k}
\\
& &{} + \int_s^{u} b'\bigl(v,
X_v^{\theta,k}\bigr) D_s X_v^{\theta,k}
\mrmd v +\int_s^{u} a'\bigl(v,
X_v^{\theta,k}\bigr) D_s X_v^{\theta,k}\mrmd W_v,
\end{eqnarray*}
and we deduce that $D_s
X_{u}^{\theta,k}=Y_{u}^{\theta,k}(1+c'(X_{T_k-}^{\theta,k},
\theta)1_{u\geq T_k})(Y_{s}^{\theta,k})^{-1}a(s, X_s^{\theta,k}) $.

It follows that:
%
\begin{equation}
\label{Eq-DX} D_s X_{(i_k+1)/{n}}^{\theta,k}=Y_{({i_k+1})/{n}}^{\theta,k}
\bigl(1+c'\bigl(X_{T_k-}^{\theta,k}, \theta\bigr)1_{s \leq T_k}\bigr)
\bigl(Y_{s}^{\theta,k} \bigr)^{-1} a\bigl(s, X_s^{\theta,k}\bigr).
\end{equation}
From (\ref{Eq-Xpoint}) and (\ref{Eq-DX}), we obtain
%
\begin{eqnarray}
\label{Eq-exPn}\hspace*{-15pt} \dot{X}_{(i_k+1)/{n}}^{\theta,k} \gamma^{\theta,k}
D_s X_{(i_k+1)/{n}}^{\theta,k}&=& \frac{ (Y_{T_k}^{\theta,k}
Y_s^{\theta,k})^{-1} (1+c'(X_{T_k-}^{\theta,k}, \theta)1_{s\leq T_k})
a(s, X_s^{\theta,k}) \dot{c}(X_{T_k-}^{\theta,k}, \theta)} {
\int_{{i_k}/{n}}^{({i_k+1})/{n}} (Y_u^{\theta,k})^{-2} a^2(u,
X_u^{\theta,k}) (1+c'(X_{T_k-}^{\theta,k}, \theta)1_{u\leq T_k})^2 \mrmd u}\nonumber\\[-4pt]\\[-12pt]
&=&
P^{n,\theta,k}_s,\nonumber
\end{eqnarray}
and the Lemma \ref{L-malliavin} is proved.
\end{pf}

In the next lemma, we explicit the conditional expectation appearing in
the decomposition of $\frac{\dot{p}^{\theta,T}}{p^{\theta,T} }
(\frac{i_k}{n}, \frac{i_k+1}{n }, x,y)$.
%
\begin{lem}\label{L-divcond}
Assuming \textup{\ref{HA1lower}} and \textup{\ref{HA2lower}}, we have
%
\begin{eqnarray}
\label{Eq-divcond}
E^{x,T,k} \bigl(\delta\bigl(P^{n,\theta,k}\bigr)
| X_{(i_k+1)/{n}}^{\theta,k}=y
\bigr)&=&\frac{(y-x-c(x,\theta))\dot{c}(x,\theta) }{D^{n,\theta,k}(x)}\nonumber\\[-8pt]\\[-8pt]
&&{} +
E^{x,T,k}\bigl(Q^{n,\theta,k} | X_{(i_k+1)/{n}}^{\theta,k}=y
\bigr)\nonumber
\end{eqnarray}
with
%
\begin{eqnarray}
\label{Eq-den}
D^{n,\theta,k}(x)&=&a^2\biggl(\frac{i_k}{n},x
\biggr) \bigl(1+c'(x,\theta)\bigr)^2
\biggl(T_k-\frac
{i_k}{n}\biggr)\nonumber\\[-8.5pt]\\[-8.5pt]
&&{}+ a^2\biggl(
\frac{i_k}{n}, x+c(x,\theta)\biggr) \biggl(\frac{i_k+1}{n}-T_k
\biggr)\nonumber
\end{eqnarray}
and where $Q^{n,\theta,k}$ satisfies
\[
\forall p \geq1\qquad\bigl(E^{x,T,k} \bigl| Q^{n, \theta,k}
\bigr|^p\bigr)^{1/p} \leq C_p
\]
for a constant $C_p$ independent of $x,n$ and $\theta$.
\end{lem}
The first term in the right-hand side of (\ref{Eq-divcond}) is the main
term and we will prove later that the contribution of the conditional
expectation of $Q^{n, \theta,k}$ is negligible.
\begin{pf}
We first give an approximation of the process $P^{n,\theta,k}$ which
depends on the position of $s$ with respect to the jump time $T_k$. We have:
%
\begin{eqnarray}
\label{Eq-delta}
P_s^{n, \theta,k}
&=&  \biggl( \bigl(1+c'\bigl(X_{{i_k}/{n}}^{\theta,k},
\theta\bigr)\bigr)
a\biggl(\frac{i_k}{n}, X_{{i_k}/{n}}^{\theta,k}\biggr) 1_{[{i_k}/{n},
T_k]}(s) \nonumber\\[-1pt]
&&\hspace*{5pt}{}+
a\biggl(\frac{i_k}{n}, X_{{i_k}/{n}}^{\theta,k}
+c\bigl(X_{i_k/n}^{\theta,k}, \theta\bigr)\biggr)1_{(T_k, ({i_k+1})/{n}]}(s) \biggr)
\dot{c}\bigl(X_{{i_k}/{n}}^{\theta,k}, \theta\bigr)\nonumber\\[-8.5pt]\\[-8.5pt]
&&{}/
D^{n,\theta,k}\bigl(X_{{i_k}/{n}}^{\theta,k}\bigr)
\nonumber\\[-1pt]
&&{} + U^{n,\theta,k}_s,\nonumber
\end{eqnarray}
where $D^{n,\theta,k}(X_{{i_k}/{n}}^{\theta,k})$ is defined by
(\ref{Eq-den}) and $ U^{n,\theta,k}_s$ is a remainder term. We deduce
then that
%
\begin{eqnarray}\label{Eq-approxdelta}
\delta\bigl(P^{n,\theta,k}\bigr)
&=& \biggl( \bigl(1+c'\bigl(X_{{i_k}/{n}}^{\theta,k},
\theta\bigr)\bigr)
a\biggl(\frac{i_k}{n}, X_{{i_k}/{n}}^{\theta,k}\biggr)(W_{T_k}-W_{{i_k}/{n}})\nonumber\\[-1pt]
&&\hspace*{5pt}{}+
a\biggl(\frac{i_k}{n}, X_{{i_k}/{n}}^{\theta,k}+c\bigl(X_{i_k/n}^{\theta,k}, \theta\bigr)\biggr)
(W_{(i_k+1)/{n}}-W_{T_k}) \biggr)
\dot{c}\bigl(X_{{i_k}/{n}}^{\theta,k}, \theta\bigr)\nonumber\\[-8.5pt]\\[-8.5pt]
&&{}/ D^{n,\theta,k}\bigl(X_{{i_k}/{n}}^{\theta,k}\bigr)
\nonumber\\[-1pt]
&&{}
+ \delta\bigl(U^{n,\theta,k}\bigr).\nonumber
\end{eqnarray}

Now, we can approximate $X_{({i_k+1})/{n}}^{\theta,k}$ in the
following way:
\begin{eqnarray*}
X_{({i_k+1})/{n}}^{\theta,k}&=&X_{{i_k}/{n}}^{\theta,k} +c
\bigl(X_{T_k-}^{\theta,k}, \theta\bigr) + a\biggl(\frac{i_k}{n},
X_{{i_k}/{n}}^{\theta,k}\biggr) (W_{T_k}-W_{i_k/n})
\\[-1pt]
&&{}+ a\biggl(\frac{i_k}{n}, X_{{i_k}/{n}}^{\theta,k}+c \bigl(X_{i_k/n}^{\theta,k}, \theta\bigr)\biggr) (W_{(i_k+1)/{n}}-W_{T_k}) +R_1^{n,\theta,k},
\end{eqnarray*}
but observing that
\[
c\bigl(X_{T_k-}^{\theta,k}, \theta\bigr)=c\bigl(X_{i_k/n}^{\theta,k},
\theta\bigr)+c'\bigl(X_{{i_k}/{n}}^{\theta,k}, \theta
\bigr)a\biggl(\frac{i_k}{n}, X_{i_k/n}^{\theta,k}\biggr)
(W_{T_k}-W_{{i_k}/{n}})+ R_2^{n,\theta,k},
\]
we finally obtain
%
\begin{eqnarray}
\label{Eq-approxX} X_{({i_k+1})/{n}}^{\theta,k}&=&X_{i_k/n}^{\theta,k}
+c\bigl(X_{{i_k}/{n}}^{\theta,k},\theta\bigr) + \bigl(1+c'
\bigl(X_{{i_k}/{n}}^{\theta,k}, \theta\bigr)\bigr)a\biggl(
\frac{i_k}{n}, X_{i_k/n}^{\theta,k}\biggr) (W_{T_k}-W_{{i_k}/{n}})\quad
\nonumber\\[-8pt]\\[-8pt]
& &{} + a\biggl(\frac{i_k}{n}, X_{{i_k}/{n}}^{\theta,k}+c
\bigl(X_{i_k/n}^{\theta,k}, \theta\bigr)\biggr) (W_{(i_k+1)/{n}}-W_{T_k})
+R^{n,\theta,k}\nonumber
\end{eqnarray}
with $R^{n,\theta,k}=R_1^{n,\theta,k}+R_2^{n,\theta,k}$.

Putting together (\ref{Eq-approxdelta}) and (\ref{Eq-approxX}), this yields
%
\begin{eqnarray}
\label{Eq-approx} \delta\bigl(P^{n,\theta,k}\bigr)&=&\frac{(
X_{({i_k+1})/{n}}^{\theta,k}-X_{{i_k}/{n}}^{\theta,k}
-c(X_{{i_k}/{n}}^{\theta,k},\theta))\dot
{c}(X_{{i_k}/{n}}^{\theta,k}, \theta)} {
D^{n,\theta,k}(X_{{i_k}/{n}}^{\theta,k})}\nonumber\\[-8pt]\\[-8pt]
&&{}-R^{n, \theta,k} \frac
{\dot{c}(X_{{i_k}/{n}}^{\theta,k}, \theta)} {
D^{n,\theta,k}(X_{{i_k}/{n}}^{\theta,k})}+
\delta\bigl(U^{n,\theta,k}\bigr).\nonumber
\end{eqnarray}
Letting $Q^{n,\theta,k}$ be the random variable defined by
%
\begin{equation}
\label{Eq-Qn} Q^{n,\theta,k}=\delta\bigl(U^{n,\theta,k}\bigr)-R^{n,
\theta,k}
\frac{\dot
{c}(X_{{i_k}/{n}}^{\theta,k}, \theta)} {
D^{n,\theta,k}(X_{{i_k}/{n}}^{\theta,k})},
\end{equation}
where $U^{n,\theta,k}$ and $R^{n,\theta,k}$ are, respectively, defined
by (\ref{Eq-delta}) and (\ref{Eq-approxX}), we deduce
easily the first part of Lemma \ref{L-divcond}. It remains to bound
$E^{x,T,k} | Q^{n,\theta,k} |^p$, $\forall p \geq1$.

We remark that from \ref{HA1lower} and \ref{HA2lower}
%
\begin{equation}
\label{Eq-BDn} 0 \leq\frac{|\dot{c}(X_{{i_k}/{n}}^{\theta,k},
\theta)
|} {
D^{n,\theta,k}(X_{{i_k}/{n}}^{\theta,k})} \leq nC
\end{equation}
for a constant $C$ independent on $n$, $k$ and $\theta$. Moreover, we have
%
\begin{eqnarray}
\label{Eq-Euler} \Bigl(E\sup_{{i_k}/{n} \leq s \leq T_k- }\bigl|
X_s^{\theta,k}
- X_{i_k/n}^{\theta,k} \bigr|^p\Bigr)^{1/p} &\leq&
\frac{C_p}{\sqrt{n}} \quad\mbox{and} \nonumber\\[-8pt]\\[-8pt]
\Bigl(E\sup_{T_k \leq s \leq
({i_k+1})/{n} }\bigl|
X_s^{\theta,k} - X_{T_k}^{\theta,k}
\bigr|^p\Bigr)^{1/p} &\leq&\frac
{C_p}{\sqrt{n}}.\nonumber
\end{eqnarray}
So, one can easily deduce that, assuming \ref{HA1lower},
\[
\bigl(E^{x,T,k} \bigl| R^{n,\theta,k} \bigr|^p
\bigr)^{1/p} \leq C_p/n,
\]
and combining this with (\ref{Eq-BDn}), we derive
\[
E^{x,T,k} \biggl(\bigl| R^{n,\theta,k} \bigr|\frac{|\dot
{c}(X_{i_k/n}^{\theta,k}, \theta) |} {
D^{n,\theta,k}(X_{{i_k}/{n}}^{\theta,k})}
\biggr)^p \leq C_p.
\]
Turning to $\delta(U^{n,\theta,k})$, we first recall that, from the
continuity property of the divergence operator (see Nualart \cite
{Nualart}, Proposition 1.5.8, page 80), we have
%
\begin{equation}
\label{Eq-Bdiv} \bigl(E^{x,T,k} \bigl|\delta\bigl(U^{n,\theta,k}\bigr)
\bigr|^p \bigr)^{1/p} \leq C_p\bigl( \bigl\|
U^{n,\theta,k} \bigr\|_p + \bigl\| DU^{n,\theta,k} \bigr\|_p\bigr),
\end{equation}
where
%
\begin{eqnarray}
\label{Eq-NU} \bigl\| U^{n,\theta,k} \bigr\|_p^p&=&
E^{x,T,k}\biggl(\int_{i_k/n}^{(i_k+1)/n} \bigl|
U_s^{n,\theta,k} \bigr|^{2} \mrmd s\biggr)^{p/2},
\\
\label{Eq-NDU} \bigl\| DU^{n,\theta,k} \bigr\|_p^p&=&E^{x,T,k}
\biggl(\int_{i_k/n}^{(i_k+1)/n}\int_{i_k/n}^{(i_k+1)/n}
\bigl| D_vU_s^{n,\theta,k} \bigr|^{2} \mrmd s\mrmd v
\biggr)^{p/2}.
\end{eqnarray}
To bound $U^{n,\theta,k}$, we first observe that from (\ref{Eq-NU})
\[
\bigl\| U^{n,\theta,k} \bigr\|_p^p \leq\biggl(
\frac{1}{n} \biggr)^{p/2} E^{x,T,k} \sup
_{{i_k}/{n} \leq s \leq({i_k+1})/{n} }\bigl| U_s^{n,\theta,k}
\bigr|^p,
\]
so we just have to prove
%
\begin{equation}
\label{EqBNUint} \Bigl(E^{x,T,k} \sup_{{i_k}/{n} \leq s \leq
({i_k+1})/{n} }\bigl|
U_s^{n,\theta,k} \bigr|^p\Bigr)^{1/p} \leq
C_p \sqrt{n}.
\end{equation}
The error term $U^{n,\theta,k}$ is defined by (\ref{Eq-delta}) as the
difference between $P_s^{\theta,n,k}$, given in (\ref{Eq-exPn}), and an
explicit ratio:
\begin{eqnarray*}
U_s^{n,\theta,k}&=& \frac{ (1+c'(X_{T_k-}^{\theta,k}, \theta)1_{s\leq
T_k}) a(s, X_s^{\theta,k}) \dot{c}(X_{T_k-}^{\theta,k}, \theta)} {
Y_{T_k}^{\theta,k} Y_s^{\theta,k}\int_{i_k/n}^{(i_k+1)/n}
(Y_u^{\theta,k})^{-2} a^2(u, X_u^{\theta,k}) (1+c'(X_{T_k-}^{\theta,k},
\theta)1_{u\leq T_k})^2 \mrmd u}
\\
&&{}- \biggl( \bigl(1+c'\bigl(X_{{i_k}/{n}}^{\theta,k}, \theta\bigr)\bigr)
a\biggl(\frac{i_k}{n}, X_{{i_k}/{n}}^{\theta,k}\biggr) 1_{[{i_k}/{n},
T_k]}(s) \\
&&\hspace*{17.7pt}{}+
a\biggl(\frac{i_k}{n}, X_{{i_k}/{n}}^{\theta,k}+c\bigl(X_{i_k/n}^{\theta,k}, \theta\bigr)\biggr)
1_{(T_k, ({i_k+1})/{n}]}(s) \biggr)
\dot{c}\bigl(X_{{i_k}/{n}}^{\theta,k}, \theta\bigr)\\
&&{}/
D^{n,\theta,k}\bigl(X_{{i_k}/{n}}^{\theta,k}\bigr).
\end{eqnarray*}
Since $\dot{c}$ and $c'$ are bounded, we see easily from (\ref
{Eq-Euler}) that the difference between the numerators is of order $1/
\sqrt{n}$.
Now,
we remark that
%
\begin{equation}
\label{Eq-Bexp} \Bigl(E\sup_{{i_k}/{n} \leq s,u \leq({i_k+1})/{n}
}\bigl| Y_{T_k}^{\theta,k}
Y_s^{\theta,k} \bigl(Y_u^{\theta,k}
\bigr)^{-2} -1 \bigr|^p\Bigr)^{1/p} \leq
\frac{C_p}{\sqrt{n}},
\end{equation}
and that, using the non-degeneracy assumption \ref{HA2lower}
%
\begin{eqnarray}
\label{EqBdeno}
&&\int_{i_k/n}^{(i_k+1)/n}
\bigl(Y_u^{\theta,k}\bigr)^{-2} a^2
\bigl(u, X_u^{\theta,k}\bigr) \bigl(1+c'
\bigl(X_{T_k-}^{\theta,k}, \theta\bigr)1_{u\leq T_k}
\bigr)^2 \mrmd u \nonumber\\[-8pt]\\[-8pt]
&&\quad\geq\frac
{\underline{a}{}^2 \min(1, \underline{a}{}^2)}{n
\sup_{{i_k}/{n} \leq u \leq({i_k+1})/{n} }
(Y_u^{\theta,k})^{2}}.\nonumber
\end{eqnarray}
So, combining (\ref{Eq-BY}), (\ref{Eq-BDn}), (\ref{Eq-Bexp}) and
(\ref{EqBdeno}), we obtain
\begin{eqnarray*}
&&\biggl(E\sup_{ s }\biggl|\frac{1}{ Y_{T_k}^{\theta,k} Y_s^{\theta,k}\int
_{i_k/n}^{(i_k+1)/n} (Y_u^{\theta,k})^{-2} a^2(u, X_u^{\theta,k})
(1+c'(X_{T_k-}^{\theta,k}, \theta)1_{u\leq T_k})^2 \mrmd u}\\
&&\hspace*{32.3pt}{} -
\frac{1}{D^{n,\theta,k}(X_{{i_k}/{n}}^{\theta,k})} \biggr|^p\biggr
)^{1/p}
\\
&&\quad\leq C_p \sqrt{n}.
\end{eqnarray*}
This proves (\ref{EqBNUint})
and consequently
%
\begin{equation}
\label{Eq-BNU} \bigl\| U^{n,\theta,k} \bigr\|_p \leq
C_p.
\end{equation}
It remains to bound the Malliavin derivative of $U^{n,\theta,k}$. From
(\ref{Eq-delta}) and (\ref{Eq-exPn}), we have
for $v \in[i_k/n, (i_k+1)/n]$
\[
D_vU_s^{n,\theta,k}= D_v P_s^{n,\theta,k}
=D_v\bigl(\dot{X}_{(i_k+1)/{n}}^{\theta,k} \gamma^{\theta,k}
D_sX_{(i_k+1)/{n}}^{\theta,k} \bigr).
\]
Under \ref{HA1lower}, the Malliavin derivatives of $\dot{X}_{(i_k+1)/{n}}^{\theta,k}$ and $D_sX_{(i_k+1)/{n}}^{\theta,k} $
are bounded
in $L^p$. Turning to the inverse of the Malliavin variance--covariance matrix
$\gamma^{\theta,k}$, given by (\ref{Eq-gamma}), we have
\[
\gamma^{\theta,k}= \frac{1}{\int_{i_k/n}^{(i_k+1)/n}
(Y^{\theta,k}_{({i_k+1})/{n}})^2(Y_u^{\theta,k})^{-2} a^2(u,
X_u^{\theta,k}) (1+c'(X_{T_k-}^{\theta,k}, \theta)1_{u\leq T_k})^2 \mrmd u}
\]
and from (\ref{Eq-BY}) and (\ref{EqBdeno}), it is easy to see that
%
\begin{eqnarray}
\label{Eq-Bgamma} \bigl(E^{x,T,k} \bigl|\gamma^{\theta,k}
\bigr|^p\bigr)^{1/p} &\leq& n C_p \quad\mbox{and}\nonumber\\[-8pt]\\[-8pt]
\Bigl(E^{x,T,k} \sup_{{i_k}/{n} \leq v \leq({i_k+1})/{n} }
\bigl| D_v
\gamma^{\theta,k} \bigr|^p\Bigr)^{1/p} &\leq& n
C_p.\nonumber
\end{eqnarray}
Putting this together, we obtain
\[
\Bigl(E^{x,T,k} \sup_{{i_k}/{n} \leq s,v \leq({i_k+1})/{n} }\bigl| D_v
U_s^{n,\theta,k} \bigr|^p\Bigr)^{1/p} \leq n
C_p
\]
and then
%
\begin{equation}
\label{Eq-BNDU} \bigl\| DU^{n,\theta,k} \bigr\|_p \leq
C_p.
\end{equation}
From (\ref{Eq-Bdiv}), (\ref{Eq-BNU}) and (\ref{Eq-BNDU}), we deduce
\[
\bigl(E^{x,T,k} \bigl|\delta\bigl(U^{n,\theta,k}\bigr)
\bigr|^p\bigr)^{1/p} \leq C_p,
\]
and the Lemma \ref{L-divcond} is proved.
\end{pf}

The bound on $Q^{n,\theta,k}$ given in Lemma \ref{L-divcond} is not
sufficient, since to obtain the LAMN property, we have to compute
the conditional expectation with $x=X_{{i_k}/{n}}^{\lambda}$ and
$y=X_{({i_k+1})/{n}}^{\lambda}$. So we complete the Lemma
\ref{L-divcond} with the following bound.

\begin{lem}\label{L-reste}
With the assumptions and notations of Lemma \ref{L-divcond}, we have
for $\theta$ such that $|\theta-\lambda_k |\leq C/\sqrt{n}$
\[
E^{x,T,k} \bigl| E^{x,T,k} \bigl(Q^{n,\theta,k} |
X_{(i_k+1)/{n}}^{\theta,k}=X_{({i_k+1})/{n}}^{\lambda} \bigr) \bigr|
\leq C',
\]
where the constant $C'$ is independent of $x,n$ and $\theta$.
\end{lem}
\begin{pf}
We first remark that
%
\begin{eqnarray}
\label{Eq-Q1}
&&E^{x,T,k} \bigl| E^{x,T,k} \bigl(Q^{n,\theta,k} |
X_{(i_k+1)/{n}}^{\theta,k}=X_{({i_k+1})/{n}}^{\lambda} \bigr)
\bigr|\nonumber\\[-8pt]\\[-8pt]
&&\quad\leq E^{x,T,k} \bigl| Q^{n,\theta,k} \bigr|\frac{p^{\lambda
_k,T}}{p^{\theta,T}}\biggl( \frac{i_k}{n}, \frac{i_k+1}{n},x,
X^{\theta,k}_{({i_k+1})/{n}}\biggr).\nonumber\vadjust{\goodbreak}
\end{eqnarray}
From H\"{o}lder's inequality and Lemma \ref{L-divcond}, we obtain for
$p> 1$, $q > 1$ such that \mbox{$1/p+1/q=1$},
%
\begin{eqnarray}
\label{Eq-Q2}
&&E^{x,T,k} \bigl| E^{x,T,k} \bigl(Q^{n,\theta,k} |
X_{(i_k+1)/{n}}^{\theta,k}=X_{({i_k+1})/{n}}^{\lambda} \bigr)
\bigr|\nonumber\\[-8pt]\\[-8pt]
&&\quad\leq C_p \biggl(E^{x,T,k}
\biggl(\frac{p^{\lambda_k,T}}{p^{\theta,T}} \biggl(\frac{i_k}{n},
\frac{i_k+1}{n},x, X^{\theta,k}_{({i_k+1})/{n}} \biggr) \biggr)^q
\biggr)^{1/q},\nonumber
\end{eqnarray}
and the result of Lemma \ref{L-reste} reduces to prove that there
exists $q_0 > 1$ such that
%
\begin{equation}
\label{Eq-QF}
E^{x,T,k}
\biggl(\frac{p^{\lambda_k,T}}{p^{\theta,T}}\biggl(\frac{i_k}{n},
\frac{i_k+1}{n},x, X^{\theta,k}_{({i_k+1})/{n}}\biggr) \biggr)^{q_0}
\leq C,
\end{equation}
where $C$ is independent of $n$, $x$ and $\theta$. We can write:
%
\begin{eqnarray}
\label{Eq-QI}
&&
E^{x,T,k} \biggl(\frac{p^{\lambda_k,T}}{p^{\theta,T}}\biggl(
\frac{i_k}{n}, \frac{i_k+1}{n},x, X^{\theta,k}_{({i_k+1})/{n}}\biggr)
\biggr)^{q_0}
\nonumber\\[-8pt]\\[-8pt]
&&\quad=
\int p^{\lambda_k,T}\biggl(\frac{i_k}{n},
\frac
{i_k+1}{n},x,y\biggr)^{q_0}p^{\theta,T}\biggl(
\frac{i_k}{n}, \frac
{i_k+1}{n},x,y\biggr)^{1-q_0}\mrmd y,\nonumber
\end{eqnarray}
and we can express the transition $p^{\theta,T}(\frac{i_k}{n}, \frac
{i_k+1}{n},x,y)$ by decomposing it in terms on the transitions of a
diffusion without jump on the time intervals $(\frac{i_k}{n}, T_k)$ and
$(T_k, \frac{i_k+1}{n})$
%
\begin{equation}
\label{Eq-trans}\hspace*{-15pt}
p^{\theta,T}\biggl(\frac{i_k}{n}, \frac{i_k+1}{n},x,y
\biggr)= \int p^{\theta,T}\biggl(\frac{i_k}{n}, T_k,x,z
\biggr)p^{\theta,T}\biggl(T_k, \frac{i_k+1}{n},z+c(z, \theta),y\biggr)
\mrmd z.
\end{equation}
Now, assuming \ref{HA1lower} and \ref{HA2lower}, we have the
following classical estimates of the transition probabilities of a
diffusion process (see Azencott \cite{Azencott84}, page 478), for some
constants $C_1$,~$C_2$:
\begin{eqnarray*}
C_1G\biggl(x,\underline{a}{}^2 \biggl(T_k- \frac{i_k}{n}\biggr), z\biggr)
&\leq& p^{\theta,T}\biggl(\frac{i_k}{n}, T_k,x,z\biggr) \leq
C_2G\biggl(x,\overline{a}{}^2 \biggl(T_k-\frac{i_k}{n}\biggr), z\biggr),
\\
C_1 G\biggl(z+c(z,\theta),\underline{a}{}^2 \biggl(
\frac{i_k+1}{n}-T_k\biggr), y\biggr) &\leq& p^{\theta,T}
\biggl(T_k, \frac{i_k+1}{n},z+c(z,\theta),y\biggr)
\\
&\leq& C_2 G\biggl(z+c(z,\theta),\overline{a}{}^2
\biggl(\frac
{i_k+1}{n}-T_k\biggr), y\biggr),
\end{eqnarray*}
where $G(m,\sigma^2,y)$ denotes the density of the Gaussian law with
mean $m$ and variance $\sigma^2$.
To simplify the notation, we note
$\sigma_{k,n}^-=T_k-\frac{i_k}{n}$ and $\sigma_{k,n}^+=\frac{i_k+1}{n}-T_k$.
Plugging this in (\ref{Eq-trans}), we obtain
%
\begin{equation}
\label{Eq-Btransinf} p^{\theta,T}\biggl(\frac{i_k}{n},
\frac{i_k+1}{n},x,y\biggr) \geq C_1 \int G\bigl(x,
\underline{a}{}^2 \sigma_{k,n}^-,z\bigr) G\bigl(z+c(z,\theta),
\underline{a}{}^2 \sigma_{k,n}^+,y\bigr) \mrmd z:=I_1.
\end{equation}
We get analogously,
%
\begin{equation}
\label{Eq-Btranssup} p^{\lambda_k,T}\biggl(\frac{i_k}{n},
\frac{i_k+1}{n},x,y\biggr) \leq C_2 \int G\bigl(x,
\overline{a}{}^2 \sigma_{k,n}^-,z\bigr) G\bigl(z+c(z,
\lambda_k), \overline{a}{}^2 \sigma_{k,n}^+,y
\bigr) \mrmd z:=I_2.
\end{equation}
Observe that, in order to bound (\ref{Eq-QI}), we have to compute an
upper bound for $p^{\lambda_k,T}$ and a lower bound for $p^{\theta,T}$,
since $1-q_0<0$.

Our aim now is to give more tractable bounds for the transition density
$p^{\theta,T}$. For this, we make the following change of variables in
the integrals $I_1$ and $I_2$ defined
in (\ref{Eq-Btranssup}) and (\ref{Eq-Btransinf}). We put $u=\varphi
(z)=z+c(z,\theta)-x-c(x,\theta)$.
We observe that $\varphi(x)=0$. Moreover, from \ref{HA1lower} and
\ref
{HA2lower}, $\varphi$ is invertible and its derivative satisfies,
for some constant $c_0$:
\[
0<\underline{a} \leq\bigl|\varphi'(z) \bigr|\leq c_0
\]
and consequently
\[
\frac{1}{c_0} | z |\leq\bigl|\varphi^{-1}(z)-
\varphi^{-1}(0) \bigr|\leq\frac{1}{\underline{a}} | z |.
\]
So we obtain, for some constant $C_1$
%
\begin{eqnarray}
\label{Eq-I1} I_1 &\geq&C_1 \int G\bigl(0,
\underline{a}{}^2 \sigma_{k,n}^-,\varphi^{-1}(u)-
\varphi^{-1}(0)\bigr) G\bigl(u+x+c(x,\theta), \underline{a}{}^2
\sigma_{k,n}^+,y\bigr) \mrmd u
\nonumber
\\
& \geq& C_1 \int G\bigl(0,\underline{a}{}^4
\sigma_{k,n}^-,u\bigr) G\bigl(x+c(x,\theta), \underline{a}{}^2
\sigma_{k,n}^+,y-u\bigr) \mrmd u
\\
& = & C_1 G\bigl(x+c(x,\theta), \underline{a}{}^4
\sigma_{k,n}^-+\underline{a}{}^2 \sigma_{k,n}^+,y
\bigr).\nonumber
\end{eqnarray}
Proceeding similarly,
%
\begin{equation}
\label{Eq-I2}
I_2 \leq C_2 G\bigl(x+c(x,
\lambda_k), c_0^2 \overline{a}{}^2
\sigma_{k,n}^-+\overline{a}{}^2 \sigma_{k,n}^+,y
\bigr).
\end{equation}
Turning back to (\ref{Eq-QI}), it follows that
%
\begin{eqnarray}
\label{Eq-QIB}
&&
E^{x,T,k} \biggl(\frac{p^{\lambda_k,T}}{p^{\theta,T}}\biggl(
\frac{i_k}{n}, \frac{i_k+1}{n},x, X^{\theta,k}_{({i_k+1})/{n}}\biggr)
\biggr)^{q_0}
\nonumber\\[-8pt]\\[-8pt]
&&\quad\leq C \int G^{ q_0}\bigl(x+c(x,\lambda_k),
\sigma_{k,n}^1,y\bigr) G^{1-q_0}\bigl(x+c(x,\theta),
\sigma_{k,n}^2,y\bigr) \mrmd y,\nonumber
\end{eqnarray}
where $\sigma_{k,n}^1=c_0^2 \overline{a}{}^2 \sigma_{k,n}^-+\overline
{a}{}^2 \sigma_{k,n}^+$ and
$\sigma_{k,n}^2=\underline{a}{}^4 \sigma_{k,n}^-+\underline{a}{}^2
\sigma
_{k,n}^+$. Since $\sigma_{k,n}^-+ \sigma_{k,n}^+=1/n$, we check that
$\sigma_{k,n}^1$ and $\sigma_{k,n}^2$ are both lower and upper bound by
some constants over $n$. Moreover, we have
\[
\sigma_{k,n}^1-\sigma_{k,n}^2=
\bigl(c_0^2 \overline{a}{}^2-\underline
{a}{}^4\bigr)\sigma_{k,n}^-+\bigl(\overline{a}{}^2-
\underline{a}{}^2\bigr)\sigma_{k,n}^+
\]
with $c_0^2 > \overline{a}{}^2$ and $\overline{a}{}^2>\underline{a}{}^2$, so
$0<\sigma_{k,n}^2/\sigma_{k,n}^1 <1$.

Turning back to the right-hand side term of (\ref{Eq-QIB}), we have to bound
\[
\int\frac{ \RMe ^{-q_0 { (y-x-c(x, \lambda_k))^2}/({2 \sigma
_{k,n}^1})}}{(2 \uppi \sigma_{k,n}^1)^{q_0/2}} \frac{ \RMe^{-(1-q_0) {
(y-x-c(x, \theta))^2}/({2 \sigma_{k,n}^2})}}{(2
\uppi \sigma_{k,n}^2)^{(1-q_0)/2}} \mrmd y
\]
with $1<q_0$. First we observe that this integral is finite if
$q_0/\sigma_{k,n}^1+(1-q_0)/\sigma_{k,n}^2 >0$, that is $ 1 < q_0 <
\sigma_{k,n}^1/(\sigma_{k,n}^1-\sigma_{k,n}^2)$. This choice of $q_0$
is possible since $0<\sigma_{k,n}^2/\sigma_{k,n}^1 <1$.
After some calculus, we get
\begin{eqnarray*}
&&
\int\frac{ \RMe^{-q_0 { (y-x-c(x, \lambda_k))^2}/({2 \sigma
_{k,n}^1})}}{(2 \uppi \sigma_{k,n}^1)^{q_0/2}} \frac{ \RMe^{-(1-q_0) {
(y-x-c(x, \theta))^2}/({2 \sigma_{k,n}^2})}}{
(2 \uppi \sigma_{k,n}^2)^{(1-q_0)/2}} \mrmd y
\\
&&\quad= \frac{\sqrt{{2 \uppi }/({ q_0/\sigma_{k,n}^1+(1-q_0)/\sigma_{k,n}^2})
} }{ (2 \uppi \sigma_{k,n}^1)^{q_0/2 }
(2 \uppi \sigma_{k,n}^2)^{(1-q_0)/2}}\RMe^{+ {c_n}(c(x,\theta)-
c(x,\lambda_k))^2/{2}}
\end{eqnarray*}
with
\[
c_n= \biggl(\frac{q_0(q_0-1)}{\sigma_{k,n}^1 \sigma_{k,n}^2} \biggr) \bigg/
\biggl(
\frac{q_0}{\sigma_{k,n}^1 }+\frac{(1-q_0)} { \sigma
_{k,n}^2} \biggr) >0.
\]
Recalling that $\sigma_{k,n}^1$ and $\sigma_{k,n}^2$ are of order
$1/n$, we observe that $c_n$ is bounded by some constant times $n$ and
assuming that $|\theta-\lambda_k |\leq C/\sqrt{n}$, we
finally obtain
\[
E^{x,T,k} \biggl(\frac{p^{\lambda_k,T}}{p^{\theta,T}}\biggl(\frac{i_k}{n},
\frac{i_k+1}{n},x, X^{\theta,k}_{({i_k+1})/{n}}\biggr) \biggr)^{q_0}
\leq C'
\]
for a constant $C'$ independent on $x$, $n$ and $\theta$ and the Lemma
\ref{L-reste} is proved.
\end{pf}

\begin{lem}\label{L-LAMN}
Assuming \textup{\ref{HA1lower}} and \textup{\ref{HA2lower}}, we have:
\begin{eqnarray*}
&& \int_{\lambda_k}^{\lambda_k+{h_k}/{\sqrt{n}}} \frac{\dot
{p}^{\theta,T}}{p^{\theta,T} } \biggl(
\frac{i_k}{n}, \frac{i_k+1}{n },X_{{i_k}/{n}}^{\lambda},
X_{(i_k+1)/{n}}^{\lambda}\biggr)\mrmd  \theta
\\
&&\quad=
h_k \frac{\sqrt{n}(X_{({i_k+1})/{n}}^{\lambda
}-X_{i_k/n}^{\lambda}-c(X_{{i_k}/{n}}^{\lambda},\lambda_k)) \dot
{c}(X_{{i_k}/{n}}^{\lambda},\lambda_k)}{nD^{n, \lambda_k,k}(
X_{{i_k}/{n}}^{\lambda})} -\frac{h_k^2}{2}
\frac{ \dot{c}(X_{{i_k}/{n}}^{\lambda
},\lambda
_k)^2}{nD^{n, \lambda_k,k}( X_{{i_k}/{n}}^{\lambda})} +\RMo_{\mathbf
{p}^{n,\lambda}}(1).
\end{eqnarray*}
\end{lem}

\begin{pf}
We deduce easily from Lemmas \ref{L-malliavin} and \ref{L-divcond} that
\begin{eqnarray*}
&&
\int_{\lambda_k}^{\lambda_k+h_k/\sqrt{n}} \frac{\dot{p}^{\theta,T}}{p^{\theta,T} } \biggl(
\frac{i_k}{n}, \frac{i_k+1}{n },x, y\biggr)\mrmd  \theta\\
&&\quad = \int
_{\lambda
_k}^{\lambda_k+h_k/\sqrt{n}} \frac{ (y-x-c(x, \theta)) \dot{c}(x,
\theta
)}{D^{n, \theta,k}( x)} \mrmd  \theta
\\
&&\qquad{} + \int_{\lambda_k}^{\lambda_k+h_k/\sqrt{n}} E^{x,T,k}
\bigl(Q^{n,\theta,k} | X_{(i_k+1)/{n}}^{\theta,k}=y\bigr) \mrmd
\theta
\end{eqnarray*}
with $(x,y)=(X_{({i_k+1})/{n}}^{\lambda},X_{i_k/n}^{\lambda})$.
From Lemma \ref{L-reste}, the second term on the right-hand side of the
preceding equation tends to zero in probability.
Now, from a Taylor expansion of $c$, we have the approximation for
$\theta\in[\lambda_k, \lambda_k+h_k/\sqrt{n}]$
%
\begin{eqnarray}
\label{Eq-epsilon}
\frac{ (y-x-c(x, \theta)) \dot{c}(x, \theta)}{D^{n,
\theta,k}( x)} & = & \frac{ (y-x-c(x, \lambda_k)-(\theta-\lambda_k)
\dot{c}(x, \lambda_k)) \dot{c}(x, \lambda_k)}{D^{n, \lambda_k,k}( x)}
\\
& &{} + \varepsilon^{n, \theta, \lambda_k}(x,y).\nonumber
\end{eqnarray}
From \ref{HA1lower}, and using (\ref{Eq-BDn}), we have $\forall
\theta\in[\lambda_k, \lambda_k+h_k/\sqrt{n}]$
%
\begin{equation}
\label{Eq-Beps1} \biggl|\frac{\dot{c}(x, \theta)}{D^{n,\theta,k}(x)} -
\frac{\dot{c}(x,
\lambda_k)}{D^{n,\lambda_k,k}(x)} \biggr|\leq C
\sqrt{n},
\end{equation}
where $C$ does not depend on $x$. So we deduce that $\forall\theta\in
[\lambda_k, \lambda_k+h_k/\sqrt{n}]$
\[
\bigl|\varepsilon^{n, \theta, \lambda_k}(x,y) \bigr|\leq C\bigl(1+ \sqrt
{n} \bigl| y-x-c(x,
\lambda_k)\bigr|\bigr)
\]
for a constant $C$ independent on $x$ and $y$.
Consequently, it follows that
\[
\int_{\lambda_k}^{\lambda_k+h_k/\sqrt{n}} \varepsilon^{n, \theta,
\lambda_k}
\bigl(X_{{i_k}/{n}}^{\lambda},X_{({i_k+1})/{n}}^{\lambda
}\bigr) \mrmd
\theta
\]
goes
to zero in probability as $n$ goes to infinity, and the thesis follows.
\end{pf}

\begin{lem}\label{L-tcl} Let us assume \textup{\ref{HA0lower}--\ref{HA2lower}}.
Let $I_n(\lambda)$ be the diagonal matrix of size $K \times K$, and
$N_n(\lambda)$ be the random vector of size $K$,
defined by the entries,
%
\begin{equation}\label{Eq-InkNnk}
I_n(\lambda)_k=\frac{ \dot{c}(X_{{i_k}/{n}}^{\lambda},\lambda
_k)^2}{nD^{n, \lambda_k,k}( X_{{i_k}/{n}}^{\lambda})},\qquad N_n(
\lambda)_k=\frac{\sqrt{n}(X_{({i_k+1})/{n}}^{\lambda
}-X_{i_k/n}^{\lambda}-c(X_{{i_k}/{n}}^{\lambda},\lambda_k))} {
\sqrt{nD^{n, \lambda_k,k}( X_{{i_k}/{n}}^{\lambda}) }}.
\end{equation}
Then, we have,
\[
\bigl(I_n(\lambda),N_n(\lambda)\bigr) \xrightarrow{
\mathit{law}} {n \to\infty} \bigl(I(\lambda),N\bigr)
\]
with $I(\lambda)$ the diagonal matrix,
\[
I(\lambda)_k=\frac{\dot{c}(X^{\lambda}_{T_k-},\lambda_k)^2}{
a^2(T_k,X^{\lambda}_{T_k-}) [1+c'(X^{\lambda}_{T_k-},\lambda_k)]^2 U_k+
a^2(T_k,X^{\lambda}_{T_k-}+c(X^{\lambda}_{T_k-},\lambda_k) )(1-U_k)},
\]
and $U=(U_1,\ldots, U_K)$ is a vector of independent uniform laws on
$[0,1]$ such that $U$, $T$ and $(W_t)_{t \in[0,1]}$ are independent,
and conditionally on $(U,T,(W_t)_{t \in[0,1]})$, $N$ is a standard
Gaussian vector in $\mathbb{R}^K$.
\end{lem}
\begin{pf}
We just have to prove the convergence in law of the couple
\[
\bigl(n D^{n, \lambda_k,k}\bigl( X_{{i_k}/{n}}^{\lambda}\bigr
),\sqrt{n}
\bigl(X_{(i_k+1)/{n}}^{\lambda}-X_{{i_k}/{n}}^{\lambda}-c
\bigl(X_{i_k/n}^{\lambda},\lambda_k\bigr)\bigr)\bigr).
\]
We have from (\ref{Eq-den})
\begin{eqnarray*}
D^{n, \lambda_k,k}\bigl( X_{{i_k}/{n}}^{\lambda}\bigr)&=& a^2
\biggl(\frac
{i_k}{n},X_{{i_k}/{n}}^{\lambda}\biggr)
\bigl(1+c'\bigl(X_{i_k/n}^{\lambda
},
\lambda_k\bigr)\bigr)^2 \biggl(T_k-
\frac{i_k}{n}\biggr)\\
&&{}+ a^2\biggl(\frac{i_k}{n}, x+c
\bigl(X_{i_k/n}^{\lambda},\lambda_k\bigr)\biggr)
\biggl(\frac{i_k+1}{n}-T_k\biggr)
\end{eqnarray*}
and from (\ref{Eq-approxX})
\begin{eqnarray*}
X_{({i_k+1})/{n}}^{\lambda}&=&X_{{i_k}/{n}}^{\lambda} +c
\bigl(X_{i_k/n}^{\lambda},\lambda_k\bigr) +
\bigl(1+c'\bigl(X_{{i_k}/{n}}^{\lambda},
\lambda_k\bigr)\bigr)a\biggl(\frac{i_k}{n}, X_{i_k/n}^{\lambda}
\biggr) (W_{T_k}-W_{{i_k}/{n}})
\\
& &{} + a\biggl(\frac{i_k}{n}, X_{{i_k}/{n}}^{\lambda}+c
\bigl(X_{i_k/n}^{\lambda}, \lambda_k\bigr)\biggr)
(W_{(i_k+1)/{n}}-W_{T_k}) +R^{n,\lambda,k},
\end{eqnarray*}
where $R^{n,\lambda,k}$ is bounded in $L^p$ by $C/n$ (see the proof of
Lemma \ref{L-divcond}). So as a straightforward consequence
of Lemma \ref{L-tcltime}, we\vadjust{\goodbreak} obtain that $(n D^{n, \lambda_k,k}(
X_{{i_k}/{n}}^{\lambda}),\sqrt{n}(X_{({i_k+1})/{n}}^{\lambda
}-X_{{i_k}/{n}}^{\lambda}-c(X_{{i_k}/{n}}^{\lambda},\lambda
_k)))$ converges in law to
\begin{eqnarray*}
&&\bigl(D^{\lambda_k,k}\bigl(X_{T_k-}^{\lambda}\bigr),
\\
&&\quad\hspace*{0pt}\bigl(1+c'\bigl(X_{T_k-}^{\lambda},
\lambda_k\bigr)\bigr)a\bigl(T_k,X_{T_k-}^{\lambda}
\bigr)\sqrt{U_k} N_k^-
+ a\bigl(T_k,
X_{T_k-}^{\lambda}+c\bigl(X_{T_k-}^{\lambda},
\lambda_k\bigr)\bigr)\sqrt{1-U_k} N_k^+
\bigr)
\end{eqnarray*}
with
\[
D^{\lambda_k,k}\bigl(X_{T_k-}^{\lambda}\bigr)=a^2
\bigl(T_k,X^{\lambda}_{T_k-}\bigr)
\bigl[1+c'\bigl(X^{\lambda}_{T_k-},
\lambda_k\bigr)\bigr]^2 U_k+
a^2\bigl(T_k,X^{\lambda}_{T_k-}+c
\bigl(X^{\lambda}_{T_k-},\lambda_k\bigr) \bigr)
(1-U_k).
\]
This gives the result of Lemma \ref{L-tcl}.
\end{pf}

As noticed earlier, the proof of Theorem \ref{Th-LAMN} follows from the
decomposition (\ref{Eq-Znk}) with $\mathbb{P}(T \in\mathcal{T}_n)
\stackrel{n \to\infty}{\hbox to 1cm{\rightarrowfill}}1$, and Lemmas \ref{L-LAMN} and
\ref{L-tcl}.

\subsection{Proof of the convolution theorem} \label{Ssconvolutionpf}
In this section, we prove the Theorem \ref{Toptimaliteconv} and some
related results.

We recall the framework described in Section~\ref{SOptim}.

$(\Omega,\mathcal{F},\mathbb{P})$ is the canonical product space,
on which are defined the independent variables $(W_t)_{t \in[0,1]}$,
$T=(T_1,\ldots,T_K)$, $\Lambda=(\Lambda_1,\ldots,\Lambda_K)$. The
probability $\mathbb{P}$ is the simple product of the corresponding
probabilities.
From this simple disintegration of the measure $\mathbb{P}$ as a
product, we can introduce $\mathbb{P}^{\lambda}$ the probability
$\mathbb{P}$ conditional on $\Lambda=\lambda\in\mathbb{R}^K$. The
process $X$ is solution of (\ref{Emodelbornesup}), and we may assume
that for any $\lambda\in\mathbb{R}^K$ the law of $X$ under $\mathbb
{P}^{\lambda}$ is
equal to the law of $X^\lambda$ solution of (\ref{Eq-par}). We recall
that $\widetilde{\Omega}$ is the extension of ${\Omega}$ which contains
the uniform variables $U_1,\ldots,U_K$, and the Gaussian variables,
$N^{-}_1,\ldots, N_K^{-}$,
$N^{+}_1,\ldots, N_K^{+}$.

With these notations, the LAMN expansion of Theorem \ref{Th-LAMN} writes,
%
\begin{eqnarray}
\label{ELAMNrecall}
&&Z_n(\lambda, \lambda+h/\sqrt{n},
T, X_{1/n},\ldots, X_{1})\nonumber\\[-8pt]\\[-8pt]
&&\quad= \sum
_{k=1}^K h_k I_n(
\lambda)_k^{1/2} N_n(\lambda)_k
- \frac{1}{2} \sum_{k=1}^K
h_k^2 I_n(\lambda)_k
+\RMo_{\mathbb
{P}^{\lambda}}(1)\nonumber
\end{eqnarray}
with
%
\begin{eqnarray}
\label{EINconv} I_n(
\lambda)_k&=&\frac{ \dot{c}(X_{{i_k}/{n}},\lambda_k)^2}{nD^{n,
\lambda_k,k}( X_{{i_k}/{n}})},\nonumber\\[-8pt]\\[-8pt]
N_n(\lambda)_k&=&\frac{\sqrt{n}(X_{({i_k+1})/{n}}-X_{i_k/n}-c(X_{{i_k}/{n}},\lambda_k))} {
\sqrt{nD^{n, \lambda_k,k}( X_{{i_k}/{n}}) }},
\nonumber\\
D^{n,\lambda_k,k}( X_{{i_k}/{n}})&=&a^2\biggl(\frac{i_k}{n},
X_{i_k/n}\biggr) \bigl(1+c'( X_{{i_k}/{n}},
\lambda_k)\bigr)^2 \biggl(T_k-
\frac{i_k}{n}\biggr)
\nonumber\\
&&{}+ a^2\biggl(\frac{i_k}{n}, X_{{i_k}/{n}}+c(
X_{{i_k}/{n}},\lambda_k)\biggr) \biggl(\frac{i_k+1}{n}-T_k
\biggr).\nonumber
\end{eqnarray}
The Theorem \ref{Th-LAMN} states the convergence in law of
$(I_n(\lambda), N_n(\lambda))$ to $(I(\lambda), N)$ under $\mathbb
{P}^\lambda$. Actually, from the proof of Lemma \ref{L-tcl}, we get the
following convergence result under $\mathbb{P}$.

\begin{prop}\label{PLAMNbis}
Assuming \textup{\ref{HA0lower}--\ref{HA2lower}}, we have the convergence
%
\begin{eqnarray}
\label{Econvbase}
&&\bigl((nT_k-i_k)_{k},
\bigl(\sqrt{n}(W_{T_k}-W_{i_k/n})\bigr)_{k}, \bigl(
\sqrt{n}(W_{(i_k+1)/n}-W_{T_k})\bigr)_{k},
I_n(\Lambda), N_n(\Lambda) \bigr)\quad
\nonumber\\[-8pt]\\[-8pt]
&&\quad\xrightarrow{\mathit{law}} {n \to\infty} \bigl((U_k)_{k},
\bigl(\sqrt{U_k} N_k^{-}
\bigr)_{k}, \bigl(\sqrt{1-U_k} N_k^{+}
\bigr)_{k}, I(\Lambda), N(\Lambda)\bigr),
\nonumber
\end{eqnarray}
where $N(\Lambda)$ is distributed as a standard Gaussian variable in
$\mathbb{R}^K$.
Moreover this convergence is stable with respect to $\mathcal{F}$,
and the last two limit variables can be represented on the extended space
$\widetilde{\Omega}$ as,
%
\begin{eqnarray}
\label{EILambda} I(\Lambda)_k&=&\frac{\dot{c}(X_{T_k-},\Lambda_k)^2}{
a^2(T_k,X_{T_k-})
(1+c'(X_{T_k-},\Lambda_k))^2 U_k+
a^2(T_k,X_{T_k})(1-U_k)},
\\
\label{ENLambda}
N(\Lambda)_k&=&\frac{a(T_k,X_{T_k-})(1+c'(X_{T_k-}, \Lambda_k))\sqrt
{U_k} N_k^-+
a(T_k, X_{T_k})\sqrt{1-U_k} N_k^+}{[a^2(T_k,X_{T_k-})
(1+c'(X_{T_k-},\Lambda_k))^2 U_k+
a^2(T_k,X_{T_k})(1-U_k)
]^{1/2}}.
\end{eqnarray}
\end{prop}

Remark that the matrix $I(\Lambda)$ is not equal to the matrix $I^\opt$
appearing in the statement of the convolution Theorem \ref
{Toptimaliteconv}. Comparing the expression (\ref{EdefI})
of $I^\opt$ with the expression (\ref{EIpara}) of $I(\lambda)$, we see
that in the parametric case, the information is proportional to $(\dot
{c}(X_{T_k-},\lambda_k))^2$.
This is quite natural.
If instead of estimating the ``mark'' $\lambda_k$ we estimate the jump,
equal to $c(X_{T_k-},\lambda_k)$ in the parametric model, we can expect
that the effect of $(\dot{c}(X_{T_k-},\lambda_k))^2$ vanishes (by a
simple first order expansion of the error of estimation). This gives
some insight on why $\dot{c}(X_{T_k-},\Lambda_k)^2$ disappears in the
expression of $I^\opt$.\looseness=-1

On the other hand, it is not immediate why the expression of the
parametric information involves the quantity $c'(X_{T_k-}, \lambda_k)$,
which is not present in the expression of $I^\opt$.
We will see that it is due to the fact that the value of the jump
$c(X_{T_k-}, \lambda_k)$ depends on the unobserved quantity $X_{T_k-}$
and thus is not a simple functional of the parameter $\lambda_k$.

If $c$ does not depend on $X$, the situation is simpler and the proof
of Theorem \ref{Toptimaliteconv}
is much easier. For this reason, in the next section we prove the
convolution theorem in this easier setting. The general proof is given
in Section~\ref{SssProofgeneral}
and some intermediate results are stated in Section~\ref{Sssintermediate}.

\subsubsection{Proof of Theorem \texorpdfstring{\protect\ref{Toptimaliteconv}}{2.1} when
\texorpdfstring{$c(x,\theta)=c(\theta)$}{c(x,theta)=c(theta)}}

We start with a simple lemma.
%
\begin{lem} \label{LstabIN}
Assume \textup{\ref{HA0lower}--\ref{HA2lower}} then for all $\lambda,h
\in
\mathbb{R}^K$,
\begin{eqnarray*}
I_n\biggl(\lambda+\frac{h}{\sqrt{n}}\biggr) -I_n(
\lambda) &\displaystyle \xrightarrow{} {n \to\infty}& 0\qquad\mbox{in $\mathbb{P}^{\lambda}$
probability,}
\\
N_n\biggl(\lambda+\frac{h}{\sqrt{n}}\biggr) -N_n(
\lambda) +I_n(\lambda)^{1/2}h &\displaystyle \xrightarrow{} {n \to\infty}&
0\qquad \mbox{in $\mathbb{P}^{\lambda}$ probability.}
\end{eqnarray*}
\end{lem}
\begin{pf}
This follows easily from the expressions (\ref{EINconv}).
\end{pf}

Assume that $\widetilde{J}{}^{n}$ is a sequence of estimators (based on
$(X_{i/n})_{i=0,\ldots,n}$) such that
\[
\sqrt{n}\bigl(\widetilde{J}{}^{n}-J\bigr)\xrightarrow{} {n \to\infty}
\widetilde{Z}
\]
in law under $\mathbb{P}$.

Then, the Theorem \ref{Toptimaliteconv} is an immediate consequence of
the following result.
%
\begin{theo}\label{TConvosimple}
Assume \textup{\ref{HA0lower}--\ref{HA3lower}} and that $c(x,\theta
)=c(\theta)$.
Denote $\dot{C}(\Lambda)$ the diagonal matrix of size $K \times K$
such that
$\dot{C}(\Lambda)_{k}=\dot{c}(\Lambda_k)$.

Then, we have the decomposition for all $n$,
%
\begin{equation}
\sqrt{n}\bigl(\widetilde{J}{}^{n}-J\bigr)=\dot{C}(\Lambda)
I_n(\Lambda)^{-1/2} N_n(\Lambda)+R_n
\end{equation}
for $(R_n)_{n}$ a sequence of random variables with values in $\mathbb{R}^K$.

Along a subsequence $(n)$ we have the convergence in law,
%
\begin{equation}
\label{Ecvsubprop} \bigl(\dot{C}(\Lambda)
I_n(\Lambda)^{-1/2}N_n(\Lambda),R_n
\bigr)\xrightarrow{} {(n) \to\infty} \bigl(\dot{C}(\Lambda) I(
\Lambda)^{-1/2}N(\Lambda),R\bigr)=\bigl(\bigl(I^\opt
\bigr)^{-1/2} N,R\bigr),
\end{equation}
where $N=N(\Lambda)$ is Gaussian, and $R$ is independent of $N$
conditionally on $I^\opt$.

In particular, we have
$\widetilde{Z}=\lim_{(n)} \sqrt{n}(\widetilde{J}{}^{n}-J)=(I^\opt
)^{-1/2} N+R$.
\end{theo}
\begin{pf}
We set $R_n=\sqrt{n}(\widetilde{J}{}^{n}-J)-\dot{C}(\Lambda)
I_n(\Lambda)^{-1/2}
N_n(\Lambda)$ and define,
%
\begin{equation}
\label{EdefWnlambda} R_n(\lambda)=
\sqrt{n}\bigl(\widetilde{J}{}^{n}-c(\lambda_k)_k
\bigr)- \dot{C}(\lambda) I_n(\lambda)^{-1/2}N_n(
\lambda),
\end{equation}
so that $R_n=R_n(\Lambda)$.
Since $\widetilde{J}{}^{n}$ is a measurable function of the $(X_{i/n})_i$,
$J=(c(\Lambda_k))_k$ and $\dot{C}(\Lambda)$ are measurable functions of
the marks, and from the expression (\ref{EINconv}),
we deduce that
$R_n=f_n((X_{i/n})_i,T,\Lambda)$ for some Borelian function $f_n$.

Using Lemma \ref{LstabIN} and the expression (\ref{EdefWnlambda}),
we easily get:
\[
R_n\biggl(\lambda+\frac{h}{\sqrt{n}}\biggr) -R_n(
\lambda) \xrightarrow{} {n \to\infty} 0\qquad \mbox{in $\mathbb{P}^{\lambda}$
probability for any $\lambda, h \in\mathbb{R}^K$}.
\]

Remark now that by Proposition \ref{PLAMNbis} and the convergence of
$\sqrt{n}(\widetilde{J}{}^{n}-J)$, we get that $(R_n)_n$ is a tight
sequence of variables.

Hence, we can apply Proposition \ref{PindepW} below. We deduce that
\[
\bigl(I_n(\Lambda), N_n(\Lambda), R_n
\bigr) \xrightarrow{\mathrm{law}} {n \to\infty} \bigl(I(\Lambda),
N(\Lambda), R
\bigr),
\]
where the limit can be represented on an extension $\widetilde{\Omega
}\times\mathbb{R}^K$ of the space $\widetilde{\Omega}$, and the
convergence is stable with respect to $(T, \Lambda, (W_t)_{t\in[0,1]})$.
On this extension, the variable $R$ is independent of $N(\Lambda)$
conditionally on $(T, \Lambda, (W_t)_{t\in[0,1]}, (U_k)_k )$. This implies
(\ref{Ecvsubprop}), and thus the theorem.
\end{pf}

\begin{prop} \label{PindepW}
Assume \textup{\ref{HA0lower}--\ref{HA3lower}}. Let
$R_n=f_n((X_{i/n})_i,T,\Lambda)) \in\mathbb{R}^K$ where $(f_n)_n$ is a
sequence of Borelian functions.
Set $R_n(\lambda)=f_n( (X_{i/n})_i,T,\lambda)$, and assume:
\begin{itemize}
\item$R_n(\lambda+\frac{h}{\sqrt{n}}) -R_n(\lambda)
\stackrel{n \to\infty}{\hbox to 1cm{\rightarrowfill}} 0$, in $\mathbb{P}^{\lambda}$
probability for any $\lambda, h \in\mathbb{R}^K$,
\item the sequence $(R_n)_n$ is tight.
\end{itemize}
Then, one has the convergence in law, along a subsequence,
%
\begin{eqnarray}
\label{EcvZindep}
&&\bigl((nT_k-i_k)_{k},
\bigl(\sqrt{n}(W_{T_k}-W_{i_k/n})\bigr)_{k}, \bigl(
\sqrt{n}(W_{(i_k+1)/n}-W_{T_k})\bigr)_{k},
I_n(\Lambda), N_n(\Lambda), R_n\bigr)\qquad
\nonumber\\[-8pt]\\[-8pt]
&&\quad\xrightarrow{\mathit{law}} {(n) \to\infty} \bigl((U_k)_{k},
\bigl(\sqrt{U_k} N_k^{-}
\bigr)_{k}, \bigl(\sqrt{1-U_k} N_k^{+}
\bigr)_{k}, I(\Lambda), N(\Lambda), R\bigr).
\nonumber
\end{eqnarray}
The limit can be represented on a extension $\widetilde{\Omega}\times
\mathbb{R}^K$ of the space $\widetilde{\Omega}$. On this space, the
variable $R$ is independent of
$N(\Lambda)$ conditionally on $(T, \Lambda, (W_t)_{t\in[0,1]},
(U_k)_k )$.
Moreover the convergence (\ref{EcvZindep}) is stable with respect to
$(T, \Lambda, (W_t)_{t\in[0,1]})$.
\end{prop}
\begin{pf}
Consider the joint law of the random variables,
%
\begin{eqnarray}
\label{Egrosvecteur}
&&\biggl(T,\Lambda,(W_t)_{t \in
[0,1]},(nT_k-i_k)_{k=1,\ldots,K},
\biggl(\frac{(W_{T_k}-W_{i_k/n})}{\sqrt{T_k-i_k/n}}\biggr)_{k=1,\ldots,K},
\nonumber\\[-8pt]\\[-8pt]
&&\quad\hspace*{0pt}\biggl(\frac{(W_{(i_k+1)/n}-W_{T_k})}{\sqrt{(i_k+1)/n-T_k}}\biggr
)_{k=1,\ldots,K}, I_n(\Lambda),
N_n(\Lambda), R_n\biggr)
\nonumber
\end{eqnarray}
defined on the corresponding canonical product space, endowed with the
usual product topology.
From the assumption, all the components of this vector are tight, and
thus the joint law is tight. Along some subsequence, it converges in
law to some limit, and thus
(\ref{EcvZindep}) holds true. The stability of the convergence with
respect to $T,\Lambda,(W_t)_{t \in[0,1]}$ is immediate. Remark that
from Proposition \ref{PLAMNbis}, the law of the limit
\[
\bigl(T,\Lambda,(W_t)_{t \in[0,1]},(U_k)_{k=1,\ldots,K},
\bigl(N_k^{-}\bigr)_{k=1,\ldots,K},
\bigl(N_k^{+}\bigr)_{k=1,\ldots,K}, I(\Lambda), N(
\Lambda), R\bigr)
\]
is known, apart for the last component $R$. It can be clearly
represented on an extension $\widetilde{\Omega}\times\mathbb{R}^K$ of
$\widetilde{\Omega}$.

To determine some information on the law of $R$, we use techniques
inspired from the proof of convolution theorems in \cite{Jeganathan82}.

Consider the following set of random variables defined on the space
$\Omega$,
%
\begin{equation}
\label{Etestfonc}
\cases{G= g(X_{s_1}, \ldots,
X_{s_r}) \qquad \mbox{with $r \ge1$ and $(s_1,\ldots,s_r) \in[0,1]^r$},
\vspace*{1pt}\cr
G_n= g(X_{{[ns_1]}/{n}}, \ldots, X_{{[ns_r]}/{n}}),
\vspace*{1pt}\cr
\kappa=k(T_1,\ldots,T_K),
\vspace*{1pt}\cr
L_n=l(nT_1-i_1,\ldots,nT_K-i_K),
\vspace*{1pt}\cr
M=m(\Lambda_1,\ldots,\Lambda_K),}
\end{equation}
where $g$, $k$, $l$, $m$ are bounded continuous functions.

For $(\mu_1,\mu_2)\in\mathbb{R}^{2K}$ we set
\[
\varphi_n(\mu_1,\mu_2)= \mathbb{E} \bigl[
\RMe^{\RMi  \mu_1 \cdot R_n} \RMe^{\RMi  \mu_2 \cdot N_n(\Lambda)} G_n \kappa L_n M
\bigr].
\]
Clearly $G_n \to G$ in probability, and from the convergence, along a
subsequence, of (\ref{Egrosvecteur}), it is simple to show
%
\begin{equation}
\label{Ecara1} \varphi_n(\mu_1,\mu_2)
\xrightarrow{} {(n) \to\infty} \mathbb{E} \bigl[ \RMe^{\RMi  \mu_1 \cdot R}
\RMe^{\RMi  \mu_2 \cdot N(\Lambda)}
G \kappa l(U_1,\ldots,U_K) M \bigr].
\end{equation}
By conditioning on the variable $\Lambda$, whose law admits a density,
we have
\[
\varphi_n(\mu_1,\mu_2)=\int
_{\mathbb{R}^K} \mathbb{E}^\lambda\bigl[ \RMe^{\RMi  \mu_1 \cdot R_n(\lambda)}
\RMe^{\RMi  \mu_2 \cdot N_n(\lambda)} G_n \kappa L_n m(\lambda) \bigr]
f_{\Lambda}(\lambda) \mrmd  \lambda.
\]
For $h \in\mathbb{R}^K$, we make a simple change of variable in the integral,
%
\begin{eqnarray*}
&&\varphi_n(\mu_1,\mu_2)
\\
&&\quad=\int_{\mathbb{R}^K} \mathbb{E}^{\lambda+h/\sqrt{n}} \bigl[
\RMe^{\RMi  \mu_1 \cdot R_n(\lambda+h/\sqrt{n})} \RMe^{\RMi  \mu_2 \cdot N_n(\lambda
+h/\sqrt{n})} G_n \kappa L_n
\bigr]m(\lambda+h/\sqrt{n}) f_{\Lambda}(\lambda+h/\sqrt{n}) \mrmd  \lambda.
\end{eqnarray*}
Now the translation is a continuous operator in $\mathbf{L}^1(\mathbb
{R})$ and by assumption
$\lambda\mapsto m(\lambda) f_\Lambda(\lambda)$ is integrable. Thus, we
easily deduce,
\[
\varphi_n(\mu_1,\mu_2)=\int
_{\mathbb{R}^K} \mathbb{E}^{\lambda
+h/\sqrt
{n}} \bigl[ \RMe^{\RMi  \mu_1 \cdot R_n(\lambda+h/\sqrt{n})}
\RMe^{\RMi  \mu_2 \cdot N_n(\lambda+h/\sqrt{n})} G_n \kappa L_n \bigr]m(\lambda)
f_{\Lambda}(\lambda) \mrmd  \lambda+ \RMo(1).
\]
From the assumptions, we know the expansion $ R_n(\lambda+h/\sqrt
{n})=R_n(\lambda)+\RMo_{\mathbb{P}^{\lambda}}(1)$, and from Lemma
\ref{LstabIN}, we have the expansion
$N_n(\lambda+\frac{h}{\sqrt{n}}) = N_n(\lambda)
-I_n(\lambda)^{1/2}h +\RMo_{\mathbb{P}^{\lambda}}(1)$.
In these expansions, all the random variables are only depending on
$((X_{i/n})_i,T)$.
But, from the LAMN property, we know that the measures $\mathbb
{P}^{\lambda}$ and $\mathbb{P}^{\lambda+h/\sqrt{n}}$, restricted to
$((X_{i/n})_i,T)$, are contiguous.
Hence, in these expansions, one can replace $\RMo_{\mathbb{P}^{\lambda
}}(1)$ with
$\RMo_{\mathbb{P}^{\lambda+h/\sqrt{n}}}(1)$. Then, using dominated
convergence theorem, one can get
\[
\varphi_n(\mu_1,\mu_2)=\int
_{\mathbb{R}^K} \mathbb{E}^{\lambda
+h/\sqrt
{n}} \bigl[ \RMe^{\RMi  \mu_1 \cdot R_n(\lambda)}
\RMe^{\RMi  \mu_2 \cdot(N_n(\lambda)-I_n(\lambda)^{1/2}h)} G_n \kappa L_n \bigr
]m(\lambda)
f_{\Lambda}(\lambda) \mrmd  \lambda+ \RMo(1).
\]
Remark that the random variables appearing in the inner expectation
only depend on the observations $((X_{i/n})_i,T)$, and thus the
likelihood ratio
$\frac{\mathbf{p}^{n,\lambda+h/\sqrt{n}}}{\mathbf{p}^{n,\lambda
}}(T,(X_{i/n})_i)=\exp(
Z_n(\lambda, \lambda+h/\sqrt{n}, T, (X_{i/n})_i)$
might be used to change the measure,
%
\begin{eqnarray}
\label{Echangproba} \varphi_n(\mu_1,
\mu_2)&=&\int_{\mathbb{R}^K} \mathbb{E}^{\lambda} \bigl[
\RMe^{\RMi  \mu_1 \cdot R_n(\lambda)} \RMe^{\RMi  \mu_2 \cdot(N_n(\lambda)-I_n(\lambda
)^{1/2}h)} \RMe^{
Z_n(\lambda, \lambda+h/\sqrt{n}, T, (X_{i/n})_i)} G_n \kappa
L_n \bigr]\quad
\nonumber\\[-8pt]\\[-8pt]
&&\hspace*{16.1pt}{}\times m(\lambda) f_{\Lambda}(\lambda) \mrmd  \lambda+ \RMo(1).
\nonumber
\end{eqnarray}
We deduce,
\[
\varphi_n(\mu_1,\mu_2)=\mathbb{E} \bigl[
\RMe^{\RMi  \mu_1 \cdot R_n} \RMe^{\RMi  \mu_2 \cdot(N_n(\Lambda)-I_n(\Lambda
)^{1/2}h)} \RMe^{
Z_n(\Lambda, \Lambda+h/\sqrt{n}, T, (X_{i/n})_i)} G_n \kappa
L_n M \bigr]+ \RMo(1).
\]
But from the LAMN expansion
(\ref{ELAMNrecall}), one can easily get
\[
Z_n\bigl(\Lambda, \Lambda+h/\sqrt{n}, T, (X_{i/n})_i
\bigr)= h^* I_n(\Lambda)^{1/2} N_n(\Lambda) -
\tfrac{1}{2} h^* I_n(\Lambda) h + \RMo_\mathbb{P}(1),
\]
where $h^*$ is the transpose of the vector $h$.
Hence, using the convergence in law of (\ref{Egrosvecteur}), and
uniform integrability of the sequence $Z_n(\Lambda, \Lambda+h/\sqrt
{n}, T, (X_{i/n})_i)$, it can be seen that
%
\begin{eqnarray}
\label{Ecara2}
&&\varphi_n(\mu_1,\mu_2)
\nonumber\\[-8pt]\\[-8pt]
&&\quad\xrightarrow{} {(n) \to\infty}
E \bigl[ \RMe^{\RMi  \mu_1 \cdot R} \RMe^{\RMi  \mu_2 \cdot(N(\Lambda)-I(\Lambda
)^{1/2}h)} \RMe^{
h^* I(\Lambda)^{1/2} N(\Lambda) - h^*
I(\Lambda) h/2} G \kappa
l(U_1,\ldots,U_K) M \bigr].\qquad\quad
\nonumber
\end{eqnarray}
Comparing the expressions (\ref{Ecara1}) and (\ref{Ecara2}), it comes
$\forall\mu_1,\mu_2,h$,
\begin{eqnarray*}
&&
E \bigl[ \RMe^{\RMi  \mu_1 \cdot R} \RMe^{\RMi  \mu_2 \cdot N(\Lambda)}
G \kappa l(U_1,\ldots,U_K) M \bigr]
\\
&&\quad=E \bigl[ \RMe^{\RMi  \mu_1 \cdot R} \RMe^{\RMi  \mu_2 \cdot(N(\Lambda)-I(\Lambda
)^{1/2}h)} \RMe^{
h^* I(\Lambda)^{1/2} N(\Lambda) - h^*
I(\Lambda) h/2} G \kappa
l(U_1,\ldots,U_K) M \bigr].
\end{eqnarray*}
We deduce that $\forall\mu_1,\mu_2,h$, the two following conditional
expectations are almost surely equal,
\begin{eqnarray*}
&&
E \bigl[ \RMe^{\RMi  \mu_1 \cdot R} \RMe^{\RMi  \mu_2 \cdot N(\Lambda)} \mid X,T, (U_k)_k,
\Lambda\bigr]
\\
&&\quad=E \bigl[ \RMe^{\RMi  \mu_1 \cdot R} \RMe^{\RMi  \mu_2 \cdot(N(\Lambda)-I(\Lambda
)^{1/2}h)} \RMe^{
h^* I(\Lambda)^{1/2} N(\Lambda) -  h^*
I(\Lambda) h/2} 
\mid X,T, (U_k)_k, \Lambda\bigr].
\end{eqnarray*}
But from continuity and analyticity arguments, it can be seen that this
equality holds, almost surely, for any $\mu_1 \in\mathbb{R}^K,\mu_2
\in\mathbb{R}^K,h \in\mathbb{C}^K$.

Hence, we can set $h=-\RMi I(\Lambda)^{-1/2} \mu_2$ in the above relation,
and find
\[
E \bigl[ \RMe^{\RMi  \mu_1 \cdot R} \RMe^{\RMi  \mu_2 \cdot N(\Lambda)} \mid X,T, (U_k)_k,
\Lambda\bigr]= E \bigl[ \RMe^{\RMi  \mu_1 \cdot R} \mid X,T, (U_k)_k,
\Lambda\bigr] \RMe^{-\mu_2^*\mu^{}_2/2}.
\]
This precisely states that, conditionally on $(X,T, (U_k)_k, \Lambda)$,
the random variables
$R$ and $N(\Lambda)$ are independent. The proposition is proved after
remarking that the Brownian motion $(W_t)_t$ can be recovered as a
measurable functional of $X,T,\Lambda$.
\end{pf}

\subsubsection{Intermediate results} \label{Sssintermediate}

The assumption $c(x,\theta)=c(\theta)$ is crucial for the proof of
Theorem \ref{TConvosimple}. Indeed if $c$ depends on the
jump-diffusion, then $J_k=c(X_{T_k-},\lambda_k)$, and instead of
(\ref{EdefWnlambda}), we have
\[
R_n(\lambda)=\sqrt{n}\bigl(\widetilde{J}{}^{n}-c(X_{T_k-},
\lambda_k)_k\bigr)- \dot{C}(X,\lambda)
I_n(\lambda)^{-1/2}N_n(\lambda),
\]
where $\dot{C}(X,\lambda)=\operatorname{diag}(\dot
{c}(X_{T_k-},\lambda
_k)_k)$. This quantity depends on $X_{T_k-}$ which is unobserved.
However, the assumption that $R_n(\lambda)$ is only function of
$((X_{i/n})_i,T)$ is essential in the Proposition \ref{PindepW} (at
the step just before equation (\ref{Echangproba})).

But if instead of $R_n(\lambda)$ we consider
\[
R^{\mathrm{obs}}_n(\lambda)=\sqrt{n}\bigl(\widetilde
{J}{}^{n}-c(X_{i_k/n},\lambda_k)_k
\bigr)- \dot{C}^\mathrm{obs}_n(\lambda) I_n(
\lambda)^{-1/2}N_n(\lambda),
\]
where
$\dot{C}^\mathrm{obs}_n(\lambda)=\operatorname{diag}(\dot
{c}(X_{i_k/n},\lambda_k)_k)$,
then, the Proposition \ref{PindepW} can be applied, and we can prove
the following modification of Theorem \ref{TConvosimple}.
%
\begin{theo}\label{Tciblediff}
Let $\widetilde{J}{}^{n}$ be any sequence of estimators such that
\[
\sqrt{n}\bigl(\widetilde{J}{}^{n}-\bigl(c(X_{i_k/n},
\Lambda_k)\bigr)_k\bigr) \xrightarrow{\mathit{law}}
{(n) \to\infty} \overline{Z}
\]
for some variable $\overline{Z}$. Then, the law of $\overline{Z}$ is
necessarily a convolution,
\[
\overline{Z}\stackrel{\mathit{law}} {=} \dot{C}(X,\Lambda) I(\Lambda
)^{-1/2} N(\Lambda)+R,
\]
where $N(\Lambda)$ is a standard Gaussian vector independent of $\dot
{C}(X,\Lambda)^{-2} I(\Lambda)$, and $R$ is some random variable
independent of $N(\Lambda)$ conditionally on $\dot{C}(X,\Lambda)^{-2}
I(\Lambda)$.
A simple expression for the entries of the diagonal matrix
$\dot{C}(X,\Lambda)^{-2} I(\Lambda)$ is
%
\begin{eqnarray}
I_{k}&=&\bigl[U_k a(T_k,X_{T_k-})^2
\bigl(1+c'(T_k,X_{T_k-})\bigr)^2\nonumber\\[-8pt]\\[-8pt]
&&\hspace*{3.5pt}{}+ (1-U_k)a(T_k,X_{T_k})^2
\bigr]^{-1}\qquad \mbox{for } k=1,\ldots,K.\nonumber
\end{eqnarray}
\end{theo}
Actually, to prove the convolution theorem when the coefficient
$c(x,\theta)$ depends on $x$, we need a strengthened version of the
Proposition \ref{PindepW}.
Indeed, we will show
that the variable $R$, in the statement of Proposition \ref{PindepW},
is independent of $N$ conditionally on any variable that can be
obtained as a limit of the observations. This yields some additional
knowledge on the dependence between the variable $R$ and the other variables.

\begin{prop} \label{PindepWrenfo}
Let us make the same assumptions as in Proposition \ref{PindepW}.
Assume furthermore that there exist a continuous function $\Psi$ with
values in $\mathbb{R}^K$
and $(A_n)_n$ a sequence of random variables depending on the
observations $(T,(X_{i/n})_i)$, such that
\begin{eqnarray*}
&&A_n - \Psi\bigl((nT_k-i_k)_{k},
\bigl(\sqrt{n}(W_{T_k}-W_{i_k/n})\bigr)_{k}, \bigl(
\sqrt{n}(W_{(i_k+1)/n}-W_{T_k})\bigr)_{k}\bigr)
\\
&&\quad\xrightarrow{} {n \to\infty} 0\qquad \mbox{in $\mathbb{P}$ probability.}
\end{eqnarray*}
Then, in the description of the limit (\ref{EcvZindep}), the
variable $R$ is independent of $N(\Lambda)$ conditionally on
$(T, \Lambda, (W_t)_{t\in[0,1]}, (U_k)_k )$ and $\Psi((U_k)_k,(\sqrt
{U_k}N_k^{-})_k,((\sqrt{1-U_k}N_k^+)_k)$.
\end{prop}
\begin{pf}
The proof is a slight modification of the proof of Proposition \ref{PindepW}.
We simply add to the list of random variables (\ref{Etestfonc}), the
new one $S_n=s(A_n)$, with $s$ being any continuous bounded function.
Accordingly, we set $\varphi_n(\mu_1,\mu_2)= \mathbb{E}
[
\RMe^{\RMi  \mu_1 \cdot R_n}
\RMe^{\RMi  \mu_2 \cdot N_n(\Lambda)}
S_n G_nKL_n M
]$. Then, the proof follows the same lines as
the proof of Proposition \ref{PindepW}.
\end{pf}

\subsubsection{Proof of Theorem \texorpdfstring{\protect\ref{Toptimaliteconv}}{2.1}. The
general case}\label{SssProofgeneral}

We prove the Theorem \ref{Toptimaliteconv} in the general situation
where $c(x,\theta)$
depends on $x$.

As seen in the previous section, a difficulty comes from the fact that
the target of the estimator
$J=(\Delta X_{T_k})_k=(c(X_{T_k-},\Lambda_k))_k$ depends on
the unobserved value $X_{T_k-}$.
We introduce $\overline{J}{}^n=(c(X_{{i_k}/{n}},\Lambda_k))_k$, and
with simple computations, one can write the following expansion, for any
sequence of estimators $\widetilde{J}{}^{n}$,
\[
\sqrt{n}\bigl( \widetilde{J}{}^{n}_k - J_k
\bigr) =\sqrt{n}\bigl( \widetilde{J}{}^{n}_k -
\overline{J}{}^n_k\bigr) - c'(X_{i_k/n},
\Lambda_k) \sqrt{n} (X_{T_k-}-X_{{i_k}/{n}}) +
\RMo_{\mathbb{P}}(1).
\]
If $\sqrt{n}( \widetilde{J}{}^{n}_k - \overline{J}{}^n_k)$ is tight we
can use
Theorem \ref{Tciblediff} and deduce,
$\lim_{(n)} \sqrt{n}( \widetilde{J}{}^{n}_k - J_k)= \widetilde{Z}_k=
\dot{C}(X,\Lambda) I(\Lambda)_k^{-1/2} N(\Lambda)_k -
c'(X_{T_k-},\Lambda_k)
a(T_k,X_{T_k-}) \sqrt{U_k} N_k^-+R_k$.
After a few algebra, involving the expressions (\ref
{EILambda})--(\ref{ENLambda}), it could be seen that this reduces to
the algebric relation
(\ref{Econvolforme}), with $N$ being some standard normal variable.
However by this method, we cannot deduce the conditional independence
of $R$ with $N$. Indeed, only the conditional independence of $R$ with
$N(\Lambda)$ is known, and we have no information about the joint law
of $R$ and $N^-$.

To solve this problem, we consider two new statistical experiments
where we add the observation of the jump-diffusion just before (or just
after) the jump. We first state the LAMN properties for these new experiments.
We omit the proof, which is similar to the proof of Theorem \ref{Th-LAMN}.

\begin{prop}[(LAMN property adding the observations before the
jumps)]\label{PLAMNsurM} Assume \textup{\ref{HA0lower}, \ref{HA1lower}}
and \textup{\ref{HA2lower}}. Denote
$(\mathbf{p}^{n,\lambda,{\mathrm{aug}^-}})$ the density on $\mathbb
{R}^{n+2K}$ of the augmented vector of observations
$\Obs^{\mathrm{aug}^-} =((X_{i/n})_i,(T_k)_k,(X_{T_k-})_k)$ under
$\mathbb{P}^\lambda$. For $\lambda\in\mathbb{R}^K$, $h \in\mathbb
{R}^K$, define the log-likelihood ratio $Z_n^{\mathrm{aug}^-}(\lambda,
\lambda+h/\sqrt{n},\Obs^{\mathrm{aug}^-}) =
\ln\frac{\mathbf{p}^{n,\lambda+h/\sqrt{n},{\mathrm
{aug}^-}}(\Obs^{\mathrm{aug}^-} )}{\mathbf
{p}^{n,\lambda,{\mathrm{aug}^-}}(\Obs^{\mathrm{aug}^-})}$.

We have the expansion:
%
\begin{eqnarray}
\label{ELAMNsurM} Z_n\bigl(\lambda, \lambda+h/
\sqrt{n}, \Obs^{\mathrm{aug}^-}\bigr)&=& \sum_{k=1}^K
h_k I_n^{\mathrm{aug}^-}(\lambda)_k^{1/2}
N_n^{\mathrm
{aug}^-}(\lambda)_k \nonumber\\[-8pt]\\[-8pt]
&&{}- \frac{1}{2}
\sum_{k=1}^K h_k^2
I^{\mathrm{aug}^-}_n(\lambda)_k +\RMo_{\mathbb
{P}^\lambda}(1),\nonumber
\end{eqnarray}
where
%
\begin{eqnarray}
\label{EINsurM}
I_n^{\mathrm{aug}^-}(\lambda)_k&=&\frac{ \dot{c}(X_{T_k-},\lambda
_k)^2}{nD^{n,
\lambda_k,k,{\mathrm{aug}^-}}( X_{T_k-})},\nonumber\\
N_n^{\mathrm{aug}^-}(
\lambda)_k&=&\frac{\sqrt{n}(X_{(i_k+1)/{n}}-X_{T_k-}-c(X_{T_k-},\lambda_k))} {
\sqrt{nD^{n, \lambda_k,k,{\mathrm{aug}^-}}( X_{T_k-}) }},
\\
D^{n,\lambda_k,k,{\mathrm{aug}^-}}( X_{T_k-})&=& a^2\bigl(T_k,
X_{T_k-}+c( X_{T_k-},\lambda_k)\bigr) \biggl(
\frac{i_k+1}{n}-T_k\biggr).\nonumber
\end{eqnarray}
Moreover,
\[
\bigl( I^{\mathrm{aug}^-}_n(\lambda), N_n^{\mathrm{aug}^-}(
\lambda)\bigr) \xrightarrow{\mathit{law}} { n \to\infty} \bigl
(I^{\mathrm{aug}^-}(
\lambda), N^{\mathrm{aug}^-}\bigr),
\]
where $I^{\mathrm{aug}^-}(\lambda)$ is the diagonal information
matrix whose
entries are
\[
I^{\mathrm{aug}^-}(\lambda)_{k}= \frac{\dot{c}(X_{T_k-},\lambda_k)^2}{
a^2(T_k, X_{T_k})(1-U_k) }
\]
and $N^{\mathrm{aug}^-}$ is a standard Gaussian vector in $\mathbb{R}^K$.
\end{prop}

\begin{prop}[(LAMN property adding the observations after the
jumps)]\label{PLAMNsurP} Assume \textup{\ref{HA0lower}, \ref{HA1lower}}
and \textup{\ref{HA2lower}}. Denote
$(\mathbf{p}^{n,\lambda,{\mathrm{aug}^+}})$ the density on $\mathbb
{R}^{n+2K}$ of the augmented vector of observations
$\Obs^{\mathrm{aug}^+} =((X_{i/n})_i,(T_k)_k,(X_{T_k})_k)$ under
$\mathbb{P}^\lambda$. For $\lambda\in\mathbb{R}^K$, $h \in\mathbb
{R}^K$, define the log-likelihood ratio $Z_n^{\mathrm{aug}^+}(\lambda,
\lambda+h/\sqrt{n},\Obs^{\mathrm{aug}^+}) =
\ln\frac{\mathbf{p}^{n,\lambda+h/\sqrt{n},{\mathrm
{aug}^+}}(\Obs^{\mathrm{aug}^+} )}{\mathbf
{p}^{n,\lambda,{\mathrm{aug}^+}}(\Obs^{\mathrm{aug}^+})}$.

We have the expansion:
%
\begin{eqnarray}
\label{ELAMNsurP} Z_n\bigl(\lambda, \lambda+h/
\sqrt{n}, \Obs^{\mathrm{aug}^+}\bigr)&=& \sum_{k=1}^K
h_k I_n^{\mathrm{aug}^+}(\lambda)_k^{1/2}
N_n^{\mathrm
{aug}^+}(\lambda)_k \nonumber\\[-8pt]\\[-8pt]
&&{}- \frac{1}{2}
\sum_{k=1}^K h_k^2
I^{\mathrm{aug}^+}_n(\lambda)_k +\RMo_{\mathbb
{P}^\lambda}(1),\nonumber
\end{eqnarray}
where
%
\begin{eqnarray}
\label{EINsurP} I_n^{\mathrm{aug}^+}(
\lambda)_k&=&\frac{ \dot{c}(X_{i_k/n},\lambda
_k)^2}{nD^{n, \lambda_k,k,{\mathrm{aug}^+}}( X_{{i_k}/{n}})},\nonumber\\
N_n^{\mathrm{aug}^+}(
\lambda)_k&=&\frac{\sqrt{n}(X_{T_k}-X_{i_k/n}-c(X_{i_k/n},\lambda_k))} {
\sqrt{nD^{n, \lambda_k,k,{\mathrm{aug}^+}}( X_{{i_k}/{n}}) }},
\\
D^{n,\lambda_k,k,{\mathrm{aug}^+}}( X_{{i_k}/{n}})&=&a^2\biggl(
\frac
{i_k}{n},X_{{i_k}/{n}}\biggr) \bigl(1+c'(X_{{i_k}/{n}},
\lambda_k)\bigr)^2 \biggl(T_k-
\frac{i_k}{n}\biggr).\nonumber
\end{eqnarray}
Moreover,
\[
\bigl( I^{\mathrm{aug}^+}_n(\lambda), N_n^{\mathrm{aug}^+}(
\lambda)\bigr) \xrightarrow{\mathit{law}} { n \to\infty} \bigl
(I^{\mathrm{aug}^+}(
\lambda), N^{\mathrm{aug}^+}\bigr),
\]
where $I^{\mathrm{aug}^+}(\lambda)$ is the diagonal information
matrix whose
entries are
\[
I^{\mathrm{aug}^+}(\lambda)_{k}= \frac{\dot{c}(X_{T_k-},\lambda_k)^2}{
a^2(T_k, X_{T_k-})(1+c'(X_{T_k-},\lambda_k))^2U_k }
\]
and $N^{\mathrm{aug}^+}$ is a standard Gaussian vector in $\mathbb{R}^K$.
\end{prop}

We now deduce convolution results from these LAMN properties.
%
\begin{prop}\label{PConvDouble}
Let $\widetilde{J}{}^{n}$ be a sequence of estimator based on the
observations of $(X_{i/n})_{i}$
and denote $\overline{J}{}^n=( c(X_{{i_k}/{n}},\Lambda))_k$. Suppose
that the sequence $\sqrt{n}(\widetilde{J}{}^{n}-\overline{J}{}^n)$ is tight
and define
$R_n^{\mathrm{aug}^-}$ and $R_n^{\mathrm{aug}^+}$ by the following expansions
%
\begin{eqnarray}
\label{EdefRsurM} \sqrt{n}\bigl(\widetilde{J}{}^{n}-
\overline{J}{}^n\bigr)
&=& \dot{C}^\mathrm{obs}_n(
\Lambda) I_n^{\mathrm{aug}^-}(\Lambda)^{-1/2}
N_n^{\mathrm{aug}^-}(\Lambda)+R_n^{\mathrm{aug}^-},
\\
\label{EdefRsurP} \sqrt{n}\bigl(\widetilde{J}{}^{n}-
\overline{J}{}^n\bigr)
&=& \dot{C}^\mathrm{obs}_n(
\Lambda) I_n^{\mathrm{aug}^+}(\Lambda)^{-1/2}
N_n^{\mathrm{aug}^+}(\Lambda)+R_n^{\mathrm{aug}^+},
\end{eqnarray}
where $I_n^{\mathrm{aug}^-}(\Lambda)$ (resp., $I_n^{\mathrm
{aug}^+}(\Lambda)$) is the diagonal
matrix with entries $(I_n^{\mathrm{aug}^-}(\Lambda)_k)$ (resp.,
$(I_n^{\mathrm{aug}^+}
(\Lambda
)_k)$) and $\dot{C}^\mathrm{obs}_n(\Lambda)$ is diagonal with
entries $\dot
{c}(X_{{i_k}/{n}},\Lambda_k)$.

Then, we have the convergence in law
%
\begin{eqnarray}
\label{Ecvlawcompli}
&&\bigl[\sqrt{n}
\bigl(X_{({i_k+1})/{n}}-X_{T_k-}-c(X_{T_k-},
\Lambda_k)\bigr)_k,\nonumber\\
&&\qquad \sqrt{n}\bigl(X_{T_k}-X_{{i_k}/{n}}-c(X_{i_k/n},
\Lambda_k)\bigr)_k, R_n^{\mathrm{aug}^-},
R_n^{\mathrm{aug}^+}\bigr]
\nonumber\\[-8pt]\\[-8pt]
&&\quad\xrightarrow{} {(n) \to\infty}\bigl[ \bigl(a(T_k,X_{T_k}) \sqrt{1-U_k}N_k^+
\bigr)_k,\nonumber\\
&&\hspace*{36.5pt}\quad \bigl(a(T_k,X_{T_k-})
\bigl(1+c'(X_{T_k-}, \Lambda_k)\bigr)
\sqrt{U_k} N_k^-\bigr)_k,
R^{\mathrm{aug}^-}, R^{\mathrm{aug}^+}\bigr].
\nonumber
\end{eqnarray}
This convergence\vspace*{1pt} holds jointly with
(\ref{Econvbase}) and the limit variables can be represented on an
extension of $\widetilde{\Omega}$.
On this space, one has, $\forall k\in\{1,\ldots,K\}$,
%
\begin{eqnarray}
\label{ErelRR} R^{\mathrm{aug}^+}_k&=&R^{\mathrm{aug}^-}_k
- a(T_k,X_{T_k-}) \bigl(1+c'(X_{T_k-},
\Lambda_k)\bigr) \sqrt{U_k} N_k^-\nonumber\\[-8pt]\\[-8pt]
&&{} +
a(T_k,X_{T_k}) \sqrt{1-U_k}
N_k^+.\nonumber
\end{eqnarray}
Moreover, conditionally on
$(T, \Lambda, (W_t)_{t\in[0,1]}, (U_k)_k, (N_k^-)_k)$, the variable
$R^{\mathrm{aug}^-}$ is independent of $(N_k^+)_k$.
In a symmetric way, conditionally on
$(T, \Lambda, (W_t)_{t\in[0,1]}, (U_k)_k, (N_k^+)_k)$, the variable
$R^{\mathrm{aug}^+}$ is independent of $(N_k^-)_k$.
\end{prop}
\begin{pf}
From the definition of the variables $R^{\mathrm{aug}^-}_n$ and
$R^{\mathrm{aug}^+}_n$ given
by equations (\ref{EdefRsurM}) and (\ref{EdefRsurP}),
we deduce immediately the relations
%
\begin{eqnarray}
\label{ERsurM} \sqrt{n}\bigl(\widetilde{J}{}^{n}-
\overline{J}{}^n\bigr) &=& \bigl[\sqrt{n}\bigl(X_{(i_k+1)/{n}}-X_{T_k-}-c(X_{T_k-},
\Lambda_k)\bigr)\bigr]_k + R_n^{\mathrm{aug}^-}+
\RMo_{\mathbb{P}}(1),
\\
\label{ERsurP} \sqrt{n}\bigl(\widetilde{J}{}^{n}-
\overline{J}{}^n\bigr) &=& \bigl[\sqrt{n}\bigl(X_{T_k}-X_{i_k/n}-c(X_{{i_k}/{n}},
\Lambda_k)\bigr)\bigr]_k+ R_n^{\mathrm{aug}^+}.
\end{eqnarray}

By a tightness argument the joint convergence, along a subsequence, of
(\ref{Econvbase}) and
(\ref{Ecvlawcompli}) is clear. The relation (\ref{ErelRR}) is a
consequence of the equality between the quantities (\ref{ERsurM}) and
(\ref{ERsurP}).

Now, we can deduce, from the LAMN property (Proposition \ref
{PLAMNsurM}), a result analogous to Proposition \ref{PindepW}. Hence
$R^{\mathrm{aug}^-}$ is independent of the limit of $N_n^{\mathrm
{aug}^-}(\Lambda)$,
conditionally on $(T, \Lambda, (W_t)_{t\in[0,1]}, (U_k)_k)$. Moreover,
remark that in the experiment $\Obs^{\mathrm{aug}^-}$, the sequence
of variables
\[
A_n=\frac{\sqrt{n} (X_{T_k-}-X_{{i_k}/{n}} ) }{a(T_k,X_{{i_k}/{n}})}
\]
is observed. But $A_n- \sqrt{n}(W_{T_k-}-W_{{i_k}/{n}} )$ converges
to zero in $\mathbb{P}$-probability. Showing a result analogous to
Proposition \ref{PindepWrenfo}, we deduce that
$R^{\mathrm{aug}^-}$ is independent of the limit of $N_n^{\mathrm
{aug}^-}(\Lambda)$,
conditionally on $(T, \Lambda, (W_t)_{t\in[0,1]}, (U_k)_k, (\sqrt{U_k}
N_k^{-})_k)$. This shows immediately that
$R^{\mathrm{aug}^-}$ is independent of $(N_k^+)$ conditionally on $(T,
\Lambda,
(W_t)_{t\in[0,1]}, (U_k)_k, (N_k^{-})_k)$, since the sigma-fields
generated by the two vectors are the same.

The conditional independence between $R^{\mathrm{aug}^+}$ and
$(N_k^-)_k$ is
obtained in a symmetric way:
one uses the LAMN property of Proposition \ref{PLAMNsurP}, and the
fact that the sequence
\[
A'_n=\frac{\sqrt{n} (X_{({i_k+1})/{n}}-X_{T_k} ) }{a(T_k,X_{T_k})}
\]
is observed in the experiment based on $\Obs^{\mathrm{aug}^+}$.
\end{pf}

Finally, we are able to prove Theorem \ref{Toptimaliteconv}.
\begin{pf*}{Proof of Theorem \ref{Toptimaliteconv}}
First, we write
%
\begin{eqnarray}\label{Epreuvefincv}
\sqrt{n}\bigl( \widetilde{J}{}^{n}_k -
J_k\bigr)&=&\sqrt{n}\bigl( \widetilde{J}{}^{n}_k
- \overline{J}{}^n_k\bigr) - \sqrt{n}\bigl(
J_k - \overline{J}{}^n_k\bigr)
\nonumber\\[-8pt]\\[-8pt]
&=&\sqrt{n}\bigl(
\widetilde{J}{}^{n}_k - \overline{J}{}^n_k
\bigr) - c'(X_{i_k/n},\Lambda_k) \sqrt{n}
(X_{T_k-}-X_{{i_k}/{n}}) + \RMo_{\mathbb{P}}(1).\nonumber
\end{eqnarray}
But the sequence $\sqrt{n}( \widetilde{J}{}^{n} - \overline{J}{}^n)$ is
tight, and we can apply Proposition
\ref{PConvDouble}. Using (\ref{ERsurM}), (\ref{Ecvlawcompli}),
and (\ref{Epreuvefincv}) we deduce
\[
\sqrt{n}\bigl(\widetilde{J}{}^{n}_k - J_k
\bigr) \xrightarrow{\mathrm{law}} {n \to\infty} - a(T_k,X_{T_k-})
c'(X_{T_k-},\Lambda_k) \sqrt{U_k}
N^{-}_k+ a(T_k,X_{T_k})
\sqrt{1-U_k} N^+_k + R^{\mathrm{aug}^-}_k.
\]
We write the last equation as
\[
\sqrt{n}\bigl(\widetilde{J}{}^{n}_k - J_k
\bigr) \xrightarrow{\mathrm{law}} {n \to\infty} a(T_k,X_{T_k-})
\sqrt{U_k} N^{-}_k+ a(T_k,X_{T_k})
\sqrt{1-U_k} N^+_k + \widetilde{R}_k,
\]
where $\widetilde{R}_k=R^{\mathrm{aug}^-}_k -( a(T_k,X_{T_k-})
(1+c'(X_{T_k-},\Lambda_k))
\sqrt{U_k} N^{-}_k )$.
Using Proposition \ref{PConvDouble}, we deduce that $\widetilde{R}$
is independent of
$N^+$ conditionally on $(T, \Lambda, (W_t)_{t\in[0,1]}, (U_k)_k, (N_k^-)_k)$.

From (\ref{ErelRR}), we have $\widetilde{R}_k=R^{\mathrm{aug}^+}_k
- (a(T_k,X_{T_k}) \sqrt{1-U_k} N^{+}_k )_k$ and we deduce that
$\widetilde
{R}$ is independent of $N^-$ conditionally on $(T, \Lambda,
(W_t)_{t\in
[0,1]}, (U_k)_k, (N_k^+)_k)$.

Remarking\vspace*{2pt} that $N^-$ and $N^+$ are independent conditionally on $(T,
\Lambda, (W_t)_{t\in[0,1]}, (U_k)_k)$, we deduce that $\widetilde{R}$
is independent of
$(N^-,N^+)$ conditionally on $(T, \Lambda, (W_t)_{t\in[0,1]}, (U_k)_k)$.

The theorem is proved.
\end{pf*}

\begin{pf*}{Proof of Corollary \ref{CKrandom}}
We introduce the conditional probability
$\widehat{\mathbb{P}}^{K_0}=\frac{\ind_{\{K=K_0 \}}} {
\mathbb{P}(K=K_0)}\mathbb{P} $ for any $K_0 \in\mathbb{N}$ such that
$\mathbb{P}(K=K_0)>0$. For any $K_0\ge0$, the sequence
$\sqrt{n}(\widetilde{J}{}^{n}-J) $ is tight (for the product topology on
$\mathbb{R}^\mathbb{N}$) under $\widehat{\mathbb{P}}^{K_0}$. So, on a
subsequence, one has the convergence in law
$\sqrt{n}(\widetilde{J}{}^{n}-J)
\displaystyle \xrightarrow{\mathrm{law}}{\widehat{\mathbb{P}}^{K_0} }
\widetilde{Z}{}^{K_0}$, moreover the subsequence may be chosen
independent of $K_0$ from a diagonal extraction argument.

Fix $K_0 \ge1$, under the probability $\widehat{\mathbb{P}}^{K_0}$,
the assumptions \ref{HA0lower}--\ref{HA3lower} are satisfied and we
can apply Theorem \ref{Toptimaliteconv} to the $K_0$ first components
of the vector $ \widetilde{Z}{}^{K_0}$. The corollary follows from the
decomposition of the law of
$\widetilde{Z}\stackrel{\mathrm{law}}{=}\sum_{K_0 \ge0} 1_{\{
K=K_0\}}
\widetilde{Z}{}^{K_0}$.\vspace*{-4pt}
\end{pf*}

\subsection{Study of the estimator \texorpdfstring{$\widehat{J}{}^{n}$}{J n}:
Proofs of Proposition \texorpdfstring{\protect\ref{Pconsistence}}{4.1} and Theorem
\texorpdfstring{\protect\ref{TTCLsaut}}{4.1}}
\label{Ssconsistencepf}\vspace*{-4pt}

\begin{pf*}{Proof of Proposition \ref{Pconsistence}}
For $k\in\{ 1, \ldots, K \}$, let us note $i_k$ the integer such that
$ i_k/n \le T_k <(i_k+1)/n$. We set $\mathcal{I}=\{ i_1,\ldots,i_K\}$
and consider a variable which counts the number of false discovery of a
jump by the estimator,
%
\begin{equation}
\label{EdefEnproofconsistence}
E_n=\sum_{i=0}^{n-1}
\ind_{\llvert X_{(i+1)/n}-X_{i/n} \rrvert\ge u_n} \ind_{ i \notin
\mathcal{I}}.
\end{equation}

For $M>0$, we define $\Omega_{M}$ as the event $\Omega_{M}=\{ \sup_{s
\in[0,1]}
[\llvert b(s,X_s)\rrvert+\llvert a(s,X_s)\rrvert] \le M \}$.

We have
%
\begin{eqnarray}\label{Eborcant1}
&&
\mathbb{P}\bigl( \{E_n \ge1\}\cap\Omega_{M}
\bigr) \nonumber\\[-1pt]
&&\quad\le \mathbb{E}[ E_n\ind_{\Omega_M} ]
\nonumber\\[-1pt]
&&\quad = \sum_{i=0}^{n-1} \mathbb{E}[
\ind_{\llvert X_{(i+1)/n}-X_{i/n}
\rrvert\ge
u_n} \ind_{ i \notin\mathcal{I}} \ind_{\Omega_M} ]
\\[-1pt]
&&\quad\le\sum_{i=0}^{n-1} \mathbb{P}
\biggl[ \biggl\{ \biggl\llvert\int_{i/n}^{({i+1})/{n}}
a(s,X_s) \mrmd  W_s + \int_{
{i}/{n}}^{({i+1})/{n}}
b(s,X_s) \mrmd s \biggr\rrvert\ge u_n \biggr\} \cap
\Omega_M \biggr]
\nonumber\\[-1pt]
&&\quad \le\sum_{i=0}^{n-1}
\mathbb{P} \biggl[ \biggl\{ \biggl\llvert\int_{i/n}^{(i+1)/n}
a(s,X_s) \mrmd  W_s \biggr\rrvert\ge u_n -
\frac{M}{n} \biggr\} \cap\Omega_M \biggr].\nonumber
\end{eqnarray}
With $a_M=(a\land M)\vee(-M)$ one has, using Markov and
Burkholder--Davis--Gundy inequalities:
%
\begin{eqnarray}\label{Eborcant2}
&&
\mathbb{P} \biggl[ \biggl\{ \biggl\llvert\int_{i/n}^{
({i+1})/{n}}
a(s,X_s) \mrmd  W_s \biggr\rrvert\ge u_n -
\frac{M}{n} \biggr\} \cap\Omega_M \biggr] \nonumber\\
&&\quad\le \mathbb{P}
\biggl[ \biggl\llvert\int_{i/n}^{({i+1})/{n}}
a_M(s,X_s) \mrmd  W_s \biggr\rrvert\ge
u_n - \frac{M}{n} \biggr]
\nonumber\\[-8pt]\\[-8pt]
&&\quad\le C_p \biggl(u_n - \frac{M}{n}
\biggr)^{-p} n^{-p/2} \qquad\forall p>0
\nonumber\\
&&\quad = C_p n^{p(\varpi-1/2)} \qquad\forall p>0.\nonumber
\end{eqnarray}
Since $\varpi<1/2$, we get, from
(\ref{Eborcant1}) and (\ref{Eborcant2}) by choosing $p$ large
enough, $\sum_{n \ge1} \mathbb{P}(
\{ E_n \ge1 \} \cap\Omega_{M})<\infty$, and by Borel\vspace*{1pt} Cantelli's
lemma we deduce that
$\mathbb{P}( \bigcap_{n\ge1} \bigcup_{p \ge n} (\{ E_p \ge1\} \cap
\Omega_{M} ) )=0$.
It immediately implies $\mathbb{P} ( (\bigcap_{n\ge1}
\bigcup_{p \ge n} \{ E_p \ge1\} ) \cap\Omega_{M} )=0$ and since
$\bigcup_{M \ge1} \Omega_{M} = \Omega$, we easily deduce that almost
surely, there exists $n$, such that $\forall p \ge n$, $E_p=0$.
Recalling the definitions (\ref{Edefihat}) and (\ref
{EdefEnproofconsistence}), we conclude that almost surely, if $n$ is
large enough,
$\{ \hat{i}_1^{n},\ldots,\hat{i}_{\widehat{K}_n}^{n} \} \subset
\mathcal
{I}$ and, as a consequence,
$\widehat{K}_n \le K$.

Now, remark that we have almost surely the convergence, for all $k \le K$,
%
\begin{equation}
\label{EinctoDelta} X_{(i_k+1)/n}-X_{i_k/n}
\xrightarrow{} {n \to\infty} X_{T_k}-X_{T_k-}=c(X_{T_k-},
\Lambda_k). 
\end{equation}
From the assumption \ref{Hpositjumpestcoef}, we have
$c(X_{T_k-},\Lambda_k)\neq0$ and using that $u_n \to0$, we deduce
that for $n$ large enough,
$\mathcal{I} \subset\{ \hat{i}_1^{n},\ldots,\hat{i}_{\widehat
{K}_n}^{n} \} $.

As a consequence, we have shown that,
%
\begin{equation}
\label{Eihatequali} \mbox{almost
surely, for $n$ large enough}\qquad \widehat{K}_n=K
\quad\mbox{and}\quad
\hat{i}^{n}_k = i_k\qquad\forall k \le K.
\end{equation}
Eventually, the proposition follows from (\ref{EdefJhat}), (\ref
{EinctoDelta}) and (\ref{Eihatequali}).
\end{pf*}

\begin{pf*}{Proof of Theorem \ref{TTCLsaut}}
We use the notation introduced in the proof of Proposition
\ref{Pconsistence}: for $k\in\{ 1, \ldots, K \}$, we have $ i_k/n \le
T_k <(i_k+1)/n$. Let us define for $1 \le k \le K$,
$G_k^{n}=X_{(i_k+1)/n}-X_{i_k/n}-\Delta X_{T_k}$ and $G_k^{n}=0$ for $k
>K$. Using (\ref{EdefJhat}) and (\ref{Eihatequali}), we see
that, almost surely, for $n$ large enough, we have $\widehat{J}{}^{n}-J=
G^{n}$. Hence, it is sufficient to study the limit in law of $\sqrt{n}
G^{n}$.

Consider any $K_0\in\mathbb{N}$ such that $\mathbb{P}(K=K_0)>0$ and
define $\widehat{\mathbb{P}}^{K_0}=\frac{\ind_{\{K=K_0 \}}} {
\mathbb
{P}(K=K_0)}\mathbb{P} $,
the conditional probability.
Actually, we will prove the convergence of $\sqrt{n} G^{n}$
conditionally on the event $\{ K=K_0 \}$, to the law of $Z$ conditional
on $\{ K=K_0 \}$, which is sufficient to prove the theorem.


For $k > K_0$ we have $G_k^{n}=0$, hence we focus only on the
components $G_k^{n}$ with $k\le K_0$.

Define $\widehat{\Omega}{}^{n}=\{\mbox{$X$ has at most one jump on each
interval of size $1/n$} \}$. We have
\[
\lim_{n \to\infty}\widehat{\mathbb{P}}^{K_0}\bigl(\widehat{
\Omega}{}^{n}\bigr)=1.
\]
On $\widehat{\Omega}{}^{n}$, the following decomposition holds true
$\widehat{\mathbb{P}}^{K_0}$ almost surely, for any $k\le K_0$,
\[
\sqrt{n}G^{n}_k=a(i_k/n,X_{i_k/n})
\alpha_{k,n}^{-}+a(T_k,X_{T_k})
\alpha_{k,n}^{+}+e_{n,k},
\]
where
%
\begin{eqnarray}
\label{Edefalpha} \alpha_{k,n}^{-}&=&
\sqrt{n} (W_{T_k}-W_{i_k/n}),\qquad\alpha_{k,n}^{+}=
\sqrt{n} (W_{(i_k+1)/n}-W_{T_k}),
\nonumber\\
e_{n,k}&=&\sqrt{n}\int_{i_k/n}^{T_k}
\bigl(a(s,X_s)-a(i_k/n,X_{i_k/n})\bigr)
\mrmd W_s
\\
&&{}+ \sqrt{n}\int^{(i_k+1)/n}_{T_k}
\bigl(a(s,X_s)-a(T_k,X_{T_k})\bigr)
\mrmd W_s + \sqrt{n}\int^{(i_k+1)/n}_{i_k/n}
b(s,X_s) \mrmd s.\nonumber
\end{eqnarray}
First, we show that $e_{n,k}$ converges to zero in $\widehat{\mathbb
{P}}^{K_0}$ probability as $n \to\infty$. Using \ref
{Hexistestcoef}, the ordinary integral converges almost surely to
zero. It remains to see that the two stochastic integrals converge to zero.

Using that the jumps times are $\mathcal{F}_0$-measurable, we can write
the stochastic integral
\[
\sqrt{n}\int_{i_k/n}^{T_k}
\bigl(a(s,X_s)-a(i_k/n,X_{i_k/n})\bigr) \mrmd W_s
\]
as a local martingale increment
\[
\int_{0}^{1} \sqrt{n} \ind_{ [i_k/n,T_k]}(s)
\bigl(a(s,X_s)-a(i_k/n,X_{i_k/n})\bigr)
\mrmd W_s.
\]
The bracket of this local martingale is
\[
\int_{i_k/n}^{T_k} n\bigl(a(s,X_s)-a(i_k/n,X_{i_k/n})
\bigr)^2 \mrmd s,
\]
which converges to zero almost surely, using the right continuity of
the process $X$.
We deduce that $\sqrt{n} \int_{i_k/n}^{T_k}
(a(s,X_s)-a(i_k/n,X_{i_k/n})) \mrmd W_s$ converge to zero in probability. We
proceed in the same way to prove that
$\sqrt{n}\int^{(i_k+1)/n}_{T_k} (a(s,X_s)-a(T_k,X_{T_k})) \mrmd W_s
\stackrel{n \to\infty}{\hbox to 1cm{\rightarrowfill}}0$ in probability. This yields to
the relation,
\[
\sqrt{n}G^{n}_k=a(i_k/n,X_{i_k/n})
\alpha_{k,n}^{-}+a(T_k,X_{T_k})
\alpha_{k,n}^{+} + \RMo_{\widehat{\mathbb{P}}^{K_0}}(1)\qquad \mbox{for $k \le
K_0$.}
\]
Using $\widetilde{\mathrm{H}0}$, and the independence between $(W_t)_{t \in
[0,1]}$ and $T$
under $\widehat{\mathbb{P}}^{K_0}$, we can apply Lem\-ma~\ref{L-tcltime}.
We get the convergence in law, under $\widehat{\mathbb{P}}^{K_0}$,
\begin{eqnarray*}
&&\bigl((T_k)_{k=1,\ldots,K_0},\bigl(\alpha_{k,n}^-
\bigr)_{k=1,\ldots,K_0},\bigl(\alpha_{k,n}^+\bigr)_{k=1,\ldots,K_0},(W_t)_{t\in[0,1]}
\bigr)
\\
&&\quad\xrightarrow{} {n \to\infty} \bigl((T_k)_{k=1,\ldots,K_0},\bigl(
\sqrt{U_k} N^{-}_k\bigr)_{k=1,\ldots,K_0},
\bigl(\sqrt{1-U_k} N^{+}_k
\bigr)_{k=1,\ldots,K_0},(W_t)_{t\in[0,1]}\bigr).
\end{eqnarray*}
Since the marks $(\Lambda_k)_k$, the Brownian motion, and the jump
times are independent,
we have that, under $\widehat{\mathbb{P}}^{K_0}$,
$ (\alpha_{k,n}^{-}, \alpha_{k,n}^{+})_{k \le K_0} $ converges in law to
$(\sqrt{U_k} N^{-}_k, \sqrt{1-U_k}N^{+}_k)_{k \le K_0}$
stably with respect to the sigma-field generated by $(W_t)_{t \in
[0,1]}$, $(T_k)_{k}$ and $(\Lambda_k)_k$. The limit can be represented
on the extended space $\widetilde{\Omega}$ endowed with the probability
$\widetilde{\mathbb{P}}$ conditional on $K=K_0$.

But the process $X$ is measurable with respect to $\mathcal{F}_1$, and
we deduce the stable convergence,
\begin{eqnarray*}
\sqrt{n}G^{n}_k&=&a(i_k/n,X_{i_k/n})
\alpha_{k,n}^{-}+a(T_k,X_{T_k})
\alpha_{k,n}^{+}
\\
&\displaystyle \xrightarrow{} {n \to\infty}& a(T_k,X_{T_k-})
\sqrt{U_k} N^{-}_k + a(T_k,X_{T_k})
\sqrt{1-U_k}N^{+}_k
\end{eqnarray*}
for $k=1,\ldots, K_0$, under $\widehat{\mathbb{P}}^{K_0}$.

By simple computations, this implies the convergence of $(\sqrt
{n}G^{n}_k)_k$ under
$\mathbb{P}$, and the theorem is proved.
\end{pf*}

\section*{Acknowledgements}

We would like to thank the referees for their careful reading and
suggestions which improved the presentation of the paper.

This research benefited partially by the support of the `Chaire Risque
de cr\'edit', F\'ed\'eration Bancaire Fran\c{c}aise.



\printhistory


\begin{thebibliography}{24}

\bibitem{AitSahalia02}
\begin{barticle}[auto:STB|2013/05/29|08:31:43]
\bauthor{\bsnm{A{\"i}t-Sahalia},~\bfnm{Yacine}\binits{Y.}}
(\byear{2002}).
\btitle{Telling from discrete data whether the underlying continuous-time model
  is a diffusion}.
\bjournal{J. Finance}
\bvolume{57}
\bpages{2075--2112}.
\bptok{imsref}%
\end{barticle}
\endbibitem

\bibitem{AitJia09}
\begin{barticle}[mr]
\bauthor{\bsnm{A{\"{\i}}t-Sahalia},~\bfnm{Yacine}\binits{Y.}},
  \bauthor{\bsnm{Fan},~\bfnm{Jianqing}\binits{J.}} \AND
  \bauthor{\bsnm{Peng},~\bfnm{Heng}\binits{H.}}
(\byear{2009}).
\btitle{Nonparametric transition-based tests for jump diffusions}.
\bjournal{J. Amer. Statist. Assoc.}
\bvolume{104}
\bpages{1102--1116}.
\bid{doi={10.1198/jasa.2009.tm08198}, issn={0162-1459}, mr={2750239}}
\bptok{imsref}%
\end{barticle}
\endbibitem

\bibitem{AitJac09AOSb}
\begin{barticle}[mr]
\bauthor{\bsnm{A{\"{\i}}t-Sahalia},~\bfnm{Yacine}\binits{Y.}} \AND
  \bauthor{\bsnm{Jacod},~\bfnm{Jean}\binits{J.}}
(\byear{2009}).
\btitle{Estimating the degree of activity of jumps in high frequency data}.
\bjournal{Ann. Statist.}
\bvolume{37}
\bpages{2202--2244}.
\bid{doi={10.1214/08-AOS640}, issn={0090-5364}, mr={2543690}}
\bptok{imsref}%
\end{barticle}
\endbibitem

\bibitem{AitJac09AOSa}
\begin{barticle}[mr]
\bauthor{\bsnm{A{\"{\i}}t-Sahalia},~\bfnm{Yacine}\binits{Y.}} \AND
  \bauthor{\bsnm{Jacod},~\bfnm{Jean}\binits{J.}}
(\byear{2009}).
\btitle{Testing for jumps in a discretely observed process}.
\bjournal{Ann. Statist.}
\bvolume{37}
\bpages{184--222}.
\bid{doi={10.1214/07-AOS568}, issn={0090-5364}, mr={2488349}}
\bptok{imsref}%
\end{barticle}
\endbibitem

\bibitem{Azencott84}
\begin{bincollection}[mr]
\bauthor{\bsnm{Azencott},~\bfnm{Robert}\binits{R.}}
(\byear{1984}).
\btitle{Densit\'e des diffusions en temps petit: D\'eveloppements
  asymptotiques. {I}}.
In \bbooktitle{Seminar on Probability, {XVIII}}.
\bseries{Lecture Notes in Math.}
\bvolume{1059}
\bpages{402--498}.
\blocation{Berlin}: \bpublisher{Springer}.
\bid{doi={10.1007/BFb0100057}, mr={0770974}}
\bptok{imsref}%
\end{bincollection}
\endbibitem

\bibitem{BarShe06}
\begin{barticle}[auto:STB|2013/05/29|08:31:43]
\bauthor{\bsnm{Barndorff-Nielsen},~\bfnm{Ole~E.}\binits{O.E.}} \AND
  \bauthor{\bsnm{Shephard},~\bfnm{Neil}\binits{N.}}
(\byear{2006}).
\btitle{Econometrics of testing for jumps in financial economics using bipower
  variation}.
\bjournal{J. Financial Econometrics}
\bvolume{4}
\bpages{1--30}.
\bptok{imsref}%
\end{barticle}
\endbibitem

\bibitem{BarSheWin06}
\begin{barticle}[mr]
\bauthor{\bsnm{Barndorff-Nielsen},~\bfnm{Ole~E.}\binits{O.E.}},
  \bauthor{\bsnm{Shephard},~\bfnm{Neil}\binits{N.}} \AND
  \bauthor{\bsnm{Winkel},~\bfnm{Matthias}\binits{M.}}
(\byear{2006}).
\btitle{Limit theorems for multipower variation in the presence of jumps}.
\bjournal{Stochastic Process. Appl.}
\bvolume{116}
\bpages{796--806}.
\bid{doi={10.1016/j.spa.2006.01.007}, issn={0304-4149}, mr={2218336}}
\bptok{imsref}%
\end{barticle}
\endbibitem

\bibitem{ConMan11}
\begin{barticle}[mr]
\bauthor{\bsnm{Cont},~\bfnm{Rama}\binits{R.}} \AND
  \bauthor{\bsnm{Mancini},~\bfnm{Cecilia}\binits{C.}}
(\byear{2011}).
\btitle{Nonparametric tests for pathwise properties of semimartingales}.
\bjournal{Bernoulli}
\bvolume{17}
\bpages{781--813}.
\bid{doi={10.3150/10-BEJ293}, issn={1350-7265}, mr={2787615}}
\bptok{imsref}%
\end{barticle}
\endbibitem

\bibitem{GloGob08}
\begin{barticle}[mr]
\bauthor{\bsnm{Gloter},~\bfnm{Arnaud}\binits{A.}} \AND
  \bauthor{\bsnm{Gobet},~\bfnm{Emmanuel}\binits{E.}}
(\byear{2008}).
\btitle{L{AMN} property for hidden processes: The case of integrated
  diffusions}.
\bjournal{Ann. Inst. Henri Poincar\'e Probab. Stat.}
\bvolume{44}
\bpages{104--128}.
\bid{doi={10.1214/07-AIHP111}, issn={0246-0203}, mr={2451573}}
\bptok{imsref}%
\end{barticle}
\endbibitem

\bibitem{Gobet01}
\begin{barticle}[mr]
\bauthor{\bsnm{Gobet},~\bfnm{Emmanuel}\binits{E.}}
(\byear{2001}).
\btitle{Local asymptotic mixed normality property for elliptic diffusion: A~{M}alliavin calculus approach}.
\bjournal{Bernoulli}
\bvolume{7}
\bpages{899--912}.
\bid{doi={10.2307/3318625}, issn={1350-7265}, mr={1873834}}
\bptok{imsref}%
\end{barticle}
\endbibitem

\bibitem{HuaTau06}
\begin{barticle}[auto:STB|2013/05/29|08:31:43]
\bauthor{\bsnm{Huang},~\bfnm{Xin}\binits{X.}} \AND
  \bauthor{\bsnm{Tauchen},~\bfnm{George}\binits{G.}}
(\byear{2006}).
\btitle{The relative contribution of jumps of total price variance}.
\bjournal{J. Financial Econometrics}
\bvolume{4}
\bpages{456--499}.
\bptok{imsref}%
\end{barticle}
\endbibitem

\bibitem{IbrHas81Book}
\begin{bbook}[mr]
\bauthor{\bsnm{Ibragimov},~\bfnm{I.~A.}\binits{I.A.}} \AND
  \bauthor{\bsnm{Has'minski{\u\i}},~\bfnm{R.~Z.}\binits{R.Z.}}
(\byear{1981}).
\btitle{Statistical Estimation: Asymptotic theory}.
\bseries{Applications of Mathematics}
\bvolume{16}.
\blocation{New York}: \bpublisher{Springer}.
\bnote{Translated from the Russian by Samuel Kotz}.
\bid{mr={0620321}}
\bptok{imsref}%
\end{bbook}
\endbibitem

\bibitem{Jacod08SPA}
\begin{barticle}[mr]
\bauthor{\bsnm{Jacod},~\bfnm{Jean}\binits{J.}}
(\byear{2008}).
\btitle{Asymptotic properties of realized power variations and related
  functionals of semimartingales}.
\bjournal{Stochastic Process. Appl.}
\bvolume{118}
\bpages{517--559}.
\bid{doi={10.1016/j.spa.2007.05.005}, issn={0304-4149}, mr={2394762}}
\bptok{imsref}%
\end{barticle}
\endbibitem

\bibitem{JacPro98}
\begin{barticle}[mr]
\bauthor{\bsnm{Jacod},~\bfnm{Jean}\binits{J.}} \AND
  \bauthor{\bsnm{Protter},~\bfnm{Philip}\binits{P.}}
(\byear{1998}).
\btitle{Asymptotic error distributions for the {E}uler method for stochastic
  differential equations}.
\bjournal{Ann. Probab.}
\bvolume{26}
\bpages{267--307}.
\bid{doi={10.1214/aop/1022855419}, issn={0091-1798}, mr={1617049}}
\bptok{imsref}%
\end{barticle}
\endbibitem

\bibitem{JacTod10}
\begin{barticle}[mr]
\bauthor{\bsnm{Jacod},~\bfnm{Jean}\binits{J.}} \AND
  \bauthor{\bsnm{Todorov},~\bfnm{Viktor}\binits{V.}}
(\byear{2010}).
\btitle{Do price and volatility jump together?}
\bjournal{Ann. Appl. Probab.}
\bvolume{20}
\bpages{1425--1469}.
\bid{doi={10.1214/09-AAP654}, issn={1050-5164}, mr={2676944}}
\bptok{imsref}%
\end{barticle}
\endbibitem

\bibitem{Jeganathan81}
\begin{barticle}[mr]
\bauthor{\bsnm{Jeganathan},~\bfnm{P.}\binits{P.}}
(\byear{1981}).
\btitle{On a decomposition of the limit distribution of a sequence of
  estimators}.
\bjournal{Sankhy\=a Ser. A}
\bvolume{43}
\bpages{26--36}.
\bid{issn={0581-572X}, mr={0656266}}
\bptok{imsref}%
\end{barticle}
\endbibitem

\bibitem{Jeganathan82}
\begin{barticle}[mr]
\bauthor{\bsnm{Jeganathan},~\bfnm{P.}\binits{P.}}
(\byear{1982}).
\btitle{On the asymptotic theory of estimation when the limit of the
  log-likelihood ratios is mixed normal}.
\bjournal{Sankhy\=a Ser. A}
\bvolume{44}
\bpages{173--212}.
\bid{issn={0581-572X}, mr={0688800}}
\bptok{imsref}%
\end{barticle}
\endbibitem

\bibitem{Kunita}
\begin{bbook}[mr]
\bauthor{\bsnm{Kunita},~\bfnm{Hiroshi}\binits{H.}}
(\byear{1997}).
\btitle{Stochastic Flows and Stochastic Differential Equations}.
\bseries{Cambridge Studies in Advanced Mathematics}
\bvolume{24}.
\blocation{Cambridge}: \bpublisher{Cambridge Univ. Press}.
\bnote{Reprint of the 1990 original}.
\bid{mr={1472487}}
\bptok{imsref}%
\end{bbook}
\endbibitem

\bibitem{Mancini04}
\begin{barticle}[mr]
\bauthor{\bsnm{Mancini},~\bfnm{Cecilia}\binits{C.}}
(\byear{2004}).
\btitle{Estimation of the characteristics of the jumps of a general
  {P}oisson-diffusion model}.
\bjournal{Scand. Actuar. J.}
\bvolume{1}
\bpages{42--52}.
\bid{doi={10.1080/034612303100170091}, issn={0346-1238}, mr={2045358}}
\bptok{imsref}%
\end{barticle}
\endbibitem

\bibitem{Mancini09}
\begin{barticle}[mr]
\bauthor{\bsnm{Mancini},~\bfnm{Cecilia}\binits{C.}}
(\byear{2009}).
\btitle{Non-parametric threshold estimation for models with stochastic
  diffusion coefficient and jumps}.
\bjournal{Scand. J. Stat.}
\bvolume{36}
\bpages{270--296}.
\bid{doi={10.1111/j.1467-9469.2008.00622.x}, issn={0303-6898}, mr={2528985}}
\bptok{imsref}%
\end{barticle}
\endbibitem

\bibitem{Nualart}
\begin{bbook}[mr]
\bauthor{\bsnm{Nualart},~\bfnm{David}\binits{D.}}
(\byear{2006}).
\btitle{The {M}alliavin Calculus and Related Topics},
\bedition{2nd} ed.
\bseries{Probability and Its Applications (New York)}.
\blocation{Berlin}: \bpublisher{Springer}.
\bid{mr={2200233}}
\bptok{imsref}%
\end{bbook}
\endbibitem

\bibitem{vandervaart}
\begin{bbook}[mr]
\bauthor{\bparticle{van~der} \bsnm{Vaart},~\bfnm{A.~W.}\binits{A.W.}}
(\byear{1998}).
\btitle{Asymptotic Statistics}.
\bseries{Cambridge Series in Statistical and Probabilistic Mathematics}
\bvolume{3}.
\blocation{Cambridge}: \bpublisher{Cambridge Univ. Press}.
\bid{mr={1652247}}
\bptok{imsref}%
\end{bbook}
\endbibitem

\bibitem{Woerner06}
\begin{bincollection}[mr]
\bauthor{\bsnm{Woerner},~\bfnm{Jeannette H.~C.}\binits{J.H.C.}}
(\byear{2006}).
\btitle{Power and multipower variation: Inference for high frequency data}.
In \bbooktitle{Stochastic Finance}
\bpages{343--364}.
\blocation{New York}: \bpublisher{Springer}.
\bid{doi={10.1007/0-387-28359-5_12}, mr={2230770}}
\bptok{imsref}%
\end{bincollection}
\endbibitem

\end{thebibliography}
\end{document}